\documentclass[11pt]{article}
\usepackage{amssymb}
\usepackage{amsthm}
\usepackage{amsfonts}
\usepackage{eucal}
\usepackage{xspace}
\usepackage{color}
\usepackage{a4wide}
\usepackage{enumerate}
\usepackage{mathdots}

\usepackage{amsmath,hhline}
\usepackage[T1]{fontenc}
\usepackage[english]{babel}
\usepackage[colorlinks=true,citecolor=red,linkcolor=blue,urlcolor=RubineRed,pdfpagetransition=Blinds,pdftoolbar=false,pdfmenubar=false]{hyperref}
\usepackage{color}
\usepackage{dsfont}
\usepackage{float}

\usepackage{tikz}
\usepackage{standalone}
\usetikzlibrary{arrows}
\usetikzlibrary{calc}
\usetikzlibrary{intersections}

\usepackage{pgfplots}

\usepackage{mathptmx}       
\usepackage{helvet}         
\usepackage{courier}        
\usepackage{type1cm}        
\usepackage{makeidx}         
\usepackage{graphicx}        
\usepackage{multicol}        
\usepackage[bottom]{footmisc}

\newtheorem{theorem}{Theorem}[section]
\newtheorem{lemma}{Lemma}[section]
\newtheorem{proposition}{Proposition}[section]
\newtheorem{corollary}{Corollary}[section]
\theoremstyle{definition}
\newtheorem{remark}{Remark}[section]

\newcommand{\field}[1]{\ensuremath{\mathbb{#1}}}
\newcommand{\C}{\field{C}\xspace}

\newcommand{\Q}{\field{Q}\xspace}
\newcommand{\R}{\field{R}\xspace}
\newcommand{\N}{\field{N}\xspace}

\DeclareMathOperator\supp{supp}

\newcommand{\MD}[1]{\textcolor{black}{#1}}

\begin{document}

\title{\bf Remarks on the controllability of parabolic systems with non-diagonalizable diffusion matrix}

\author{Michel Duprez\thanks{MIMESIS team, Inria Nancy - Grand Est and MLMS team, Universit\'e de Strasbourg, France. E-mail: \texttt{michel.duprez@inria.fr}} 
\and Manuel Gonz\'alez-Burgos\thanks{Dpto.~Ecuaciones Diferenciales y An\'alisis Num\'erico and Instituto de Matem\'aticas 	de la Universidad de Sevilla (IMUS), Facultad de Matem\'aticas, Universidad de Sevilla, C/ Tarfia S/N, 41012 Sevilla, 	Spain.  E-mail: \texttt{manoloburgos@us.es}} 
\and Diego A. Souza\thanks{Dpto.~Ecuaciones Diferenciales y An\'alisis Num\'erico and Instituto de Matem\'aticas 	de la Universidad de Sevilla (IMUS), Facultad de Matem\'aticas, Universidad de Sevilla, C/ Tarfia S/N, 41012 Sevilla, 	Spain.  E-mail: \texttt{desouza@us.es}}}
%
%
%
%
\maketitle
\abstract{The distributed null controllability for coupled parabolic systems with non-diagonalizable diffusion matrices with a reduced number of controls has been studied in the case of constant matrices. On the other hand, boundary controllability issues and distributed controllability with non-constant coefficients for 
\MD{this} kind of systems is not completely understood. In this paper, we analyze the boundary controllability properties of a class of coupled parabolic systems with non-diagonalizable diffusion matrices in the constant case and the distributed controllability of a $2\times 2$ non-diagonalizable parabolic system with space-dependent coefficients. For the boundary controllability problem, our strategy relies on the moment method. For the distributed controllability problem, our findings provide positive and negative control results by using the Fattorini-Hautus test and a fictitious control strategy.
}
	%

%


\section{Introduction and main results}

Let  $T>0$ be a fixed time and $\omega:=(a,b)$ be a non empty open subinterval of $(0,\pi)$. Hereafter, we shall use the notation $Q_T:=(0,T)\times(0,\pi)$. We denote by $M^\star$ the conjugate transpose of a matrix $M$ and by $e_i$ the $i$-th element of the canonical basis of $\mathbb{R}^n$ ($n\in\mathbb{N}^\star$, {with $n\geq 2$}, will be fixed later). On the other hand, we also denote $(\cdot, \cdot)_{\R^n}$ (resp., $(\cdot, \cdot)_{\C^n}$) the scalar product in $\R^n$ (resp., the hermitian product in $\C^n$) and $| \cdot |$ is the euclidean norm in $\R^n$ or $\C^n$.
	
	In this work, we consider the following controlled parabolic systems in which the diffusion matrix is non-diagonalizable
\begin{equation}\label{eq:syst int}
	\left\{
	\begin{array}{ll}
		y_t  -Dy_{xx} + Ay=B{1}_{\omega}u	&\hbox{ in } Q_T,\\\noalign{\smallskip}
		y(\cdot,0)=y(\cdot,\pi)= 0			&\hbox{ in } (0,T),\\\noalign{\smallskip}
		y(0,\cdot)= y^0					&\hbox{ in } (0,\pi)
	\end{array}
	\right.
\end{equation}
and
\begin{equation}\label{eq:syst bound}
	\left\{
	\begin{array}{ll}
		z_t  -Dz_{xx} + Az=0	&\hbox{in } Q_T,\\
		\noalign{\smallskip}
		z(\cdot,0)= Bv , \quad z(\cdot,\pi)= 0 & \hbox{in } (0,T),\\\noalign{\smallskip}
		z(0,\cdot)= z^0 &\hbox{in } (0,\pi),
	\end{array}
	\right.
\end{equation}
	where $y^0\in L^2(0,\pi;\R^n)$ and $z^0\in H^{-1}(0,\pi;\R^n)$ are the initial data, 
	$u\in L^2(Q_T;\R^m)$ ($m\in \mathbb{N}^\star$) is the distributed control, $v\in L^2(0,T;\R^m)$ is the boundary control, $A\in \mathcal{L}(\R^n)$ is a zero order coupling matrix, $B\in\mathcal{L}(\R^m;\R^n)$ 
	is a matrix through which the controls act on the system and $D\in \mathcal{L}(\R^n)$ is a
	non-diagonalizable diffusion matrix satisfying an ellipticity condition given by
	$$
\left( D \xi,\xi \right)_{\R^n} \geq \beta |\xi|^2,\quad \forall \xi \in \R^{n},
	$$
	with $\beta >0$.
	
	It is well-known that for any initial datum $y^0\in L^2(0,\pi;\R^n)$ (resp., $z^0\in H^{-1}(0,\pi;\R^n)$)
	and any control $u\in L^2(Q_T;\R^m)$ (resp., $v\in L^2(0,T;\R^m)$), system \eqref{eq:syst int} (resp., 
	system \eqref{eq:syst bound}) admits a unique weak solution $y$ (resp., a unique solution by transposition $z$) with the regularity
	\begin{equation*}
	\begin{array}{c}
y\in L^2(0,T;H_0^1(0,\pi;\R^n))\cap  \mathcal{C}^0([0,T]; L^2(0,\pi;\R^n)) \\
	\noalign{\smallskip}
(\hbox{resp}.~z\in L^2(Q_T;\R^n)\cap  \mathcal{C}^0([0,T]; H^{-1}(0,\pi;\R^n))).
	\end{array}
\end{equation*}
	For more details, see \cite[p. 102]{lionscontrole} and \cite[Prop. 2.2]{fernandezcaraboundary2010}.

	Let us recall the different concepts of controllability for \eqref{eq:syst int} and \eqref{eq:syst bound} that we study in the present paper:
\begin{enumerate}
\item[$\bullet$] 
	system \eqref{eq:syst int} (resp., system \eqref{eq:syst bound}) is approximately controllable at time $T$ 
	if for any $y^0$ and $y^T$ in $  L^2(0,\pi;\R^n)$ (resp., $z^0$ and $z^T$ in $ H^{-1}(0,\pi;\R^n)$) and any $\varepsilon>0$,
	there exists a control $u\in L^2(Q_T;\R^m)$ (resp., $v\in L^2(0,T;\R^m)$) such that the associated solution 
	$y$ (resp., $z$) to system \eqref{eq:syst int} (resp., system \eqref{eq:syst bound}) satisfies
	$$
	\|y(T, \cdot)-y^T\|_{L^2(0,\pi;\R^n)}\leq \varepsilon, \quad (\hbox{resp.,}~\|z(T, \cdot)-z^T\|_{H^{-1}(0,\pi;\R^n)}\leq \varepsilon).
	$$
\item[$\bullet$]
	system \eqref{eq:syst int} (resp., system \eqref{eq:syst bound}) is null controllable at time $T$ if for any $y^0\in  L^2(0,\pi;\R^n)$ 
	(resp., $z^0\in H^{-1}(0,\pi;\R^n)$), there exists a control $u\in L^2(Q_T;\R^m)$ (resp., $v\in L^2(0,T;\R^m)$)
	such that the associated solution $y$ (resp., $z$) to system \eqref{eq:syst int} (resp., system \eqref{eq:syst bound}) satisfies
	$$
	y(T, \cdot)=0~\hbox{ in }(0,\pi), \quad (\hbox{resp.,}~z(T, \cdot)=0~\hbox{ in }(0,\pi)).
	$$
\end{enumerate}

\smallskip

 	Concerning the case of systems of ordinary differential equations under the form
\begin{equation}\label{eq:syst ode}
	\left\{
	\begin{array}{lll}
		y_t   + Ay=Bu	&~\hbox{in}& (0,T),\\\noalign{\smallskip}
		y(0)= y^0,						&
	\end{array}
	\right.
\end{equation}
where $y^0\in\R^n$, $A\in \mathcal{L}(\R^n)$ and $B\in\mathcal{L}(\R^m;\R^n)$, it is well-known that the exact controllability for~\eqref{eq:syst ode} is equivalent to the so-called \textit{Kalman rank condition}
	\begin{equation}\label{kalman}
\hbox{rank}\, [A\, : \, B]=n,
	\end{equation}
where $[A\, : \, B] \in \mathcal{L}(\R^{nm} ; \R^n)$ is the matrix given by
	$$
[A\, : \, B] := [ B \,| \, AB \, | \, \cdots \, | \, A^{n-1} B ].
	$$
This result was proved in \cite{KFA69}.

Concerning systems of partial differential equations, precisely PDEs of parabolic type, the first results about null controllability of the heat equation (\eqref{eq:syst int} and \eqref{eq:syst bound} with $n=1$), have been established in the one-dimensional case through the moment method by H.O. Fattorini and D.L. Russell, see \cite{fattorini_russell}. The distributed null controllability of the heat equation in the multi-dimensional case, has been established later, simultaneously, by G. Lebeau and L. Robbiano in \cite{LR95}, using local elliptic Carlerman estimates, and by A. Fursikov and O. Yu. Imanuvilov in \cite{FI96}, using global parabolic Carleman estimates. Using an extension method, it is possible to prove that the internal null controllability and the boundary null controllability are equivalent for the heat equation and, in general, for scalar parabolic problems (see, for instance,~\cite{AKBGBdT11bis}). It is interesting to point out that, in the context of scalar parabolic partial differential equations, the controllability properties are valid for any distributed or boundary control domain and for any time $T > 0$, i.e., there is no minimal time for controllability and no geometric restrictions on the internal or boundary control domains, contrarily to the wave equation and transport equation.

To the authors' knowledge, there are not many works devoted to the controllability of coupled parabolic systems. Unlike the scalar case, in~\cite{fernandezcaraboundary2010}, it was proved that the equivalence between the controllability of the systems of parabolic equations~\eqref{eq:syst int} and~\eqref{eq:syst bound} does not hold ({details will be provided below}) and their controllability issues should be analyzed separably. {Almost all the papers in the literature are devoted to the controllability of parabolic systems with distributed controls, acting on a small open region $\omega$ of the domain $\Omega \subset \R^N$; see, for example,~\cite{GB-dT,AKBDGB09bis,AKBDGB09,FCGBdT16,lissy_zuazua}. About 
\MD{another} kind of systems, for instance, some boundary controllability results for a system of wave equations and distributed controllability results for hyperbolic systems of first-order have been obtained in \cite{fatiha,LL17,APT20} and \cite{CORON_OLIVE,Coron-Nguyen}, respectively.} 

 Let us describe the state of the art in the case of parabolic systems with diagonalizable and non-diagonalizable diffusion matrices.


\paragraph{Diagonalizable diffusion matrices:}
In~\cite{AKBDGB09}, the authors have proved, in the case of diagonalizable diffusion matrices $D$,  that system \eqref{eq:syst int} (constant coefficients and distributed controls) is null (resp., approximate) controllable at time $T$ if and only if 
	\begin{equation}\label{kalman spec}
~\hbox{rank}\, \left[ \mu_kD+A \, : \, B \right] = n~\hbox{ for all }k\in\mathbb{N^{\star}},
	\end{equation}
where $\mu_k$ are the {eigenvalues of $-\partial_{xx}$ in $(0,\pi)$} with homogeneous Dirichlet boundary conditions. When the matrix $D$ is equal to the identity, conditions~\eqref{kalman} and~\eqref{kalman spec} are equivalent. It is surprising to obtain the same condition as in finite dimension. We refer to \cite{AKBDGB09bis} for a study in the case of time dependent coupling matrices.

{
The case of coupling matrices depending on $(t, x)$ is more intricate but in some particular parabolic systems it is possible to prove a null controllability result (cascade systems, see~\cite{GB-dT}). Let us describe the existing results for system~\eqref{eq:syst int} when $n= 2$, $m = 1$, $D = \text{diag} \, (d_1, d_2)$ ($d_1, d_2 > 0$) and 
    $$
A = A(t, x) \in L^\infty (Q_T; \mathcal{L} (\R^2)) \quad \hbox{and} \quad
B=
	\left[
	\begin{array}{c}
			 0		\\
			 1
		\end{array}
	\right],
    $$
i.e., for the $2 \times 2$ system:
	\begin{equation}\label{cascade}
	\left\{\begin{array}{ll}
\partial_t y_1 - d_ 1 \partial_{xx} y_1 + a_{11}  y_1 + a_{12}  y_2 = 0 &\hbox{in } Q_T,\\
	\noalign{\smallskip}
\partial_t y_2 - d_2 \partial_{xx} y_2 + a_{21} y_1 + a_{22} y_ 2 = {1}_{\omega}u &\hbox{in } Q_T,\\
	\noalign{\smallskip}
y_1(\cdot,0)=y_2(\cdot,0)=y_1(\cdot,\pi)=y_2(\cdot,\pi)=0&\hbox{in } (0,T),\\
	\noalign{\smallskip}
y_1(0,\cdot)=y^0_1, \quad y_2(0,\cdot)=y^0_2&\hbox{in }(0,\pi),
	\end{array}\right.
	\end{equation}
where $a_{ij} \in L^{\infty}(Q_T)$ are given functions ($1 \le i, j \le 2$), $y^0 = (y^0_1,y^0_2)^\star \in L^2(0,\pi;\mathbb{R}^2)$ is the initial datum and $u \in L^2(Q_T)$ is the distributed control. Under the generic assumption
    \begin{equation}\label{cascade-cond}
a_{12}\ge a_0>0  \mbox{ or } {-a_{12}}\ge a_0>0 \mbox{ in }  (0,T) \times \omega_0,\quad \forall i: 2\le i \le m.
    \end{equation}
for an open set $\omega_0 \subset \omega$, in~\cite{GB-dT} the authors prove a null controllability result at time $T $ for system~\eqref{cascade} which is independent of $d_1$, $d_2$, $a_{11}$, $a_{21}$ and $a_{22}$. 

The situation strongly changes if $\text{Supp}\, a_{12} \cap \omega = \emptyset $ (see for instance~\cite{boyer_olive}, \cite{CherifMinimalTimeDisjoint} and~\cite{CherifMinimalTimeDisjointB}) or the coupling term $a_{12} y_2 $ is changed by a first order coupling term (see~\cite{O14,Duprezlissy1,Duprezlissy2}). In the first case, system~\eqref{cascade} could have a minimal controllability time $T_0 \in [0, +\infty] $ such that if $T < T_0$ the system is not null controllable at time $T$ and it is if $T > T_0$ (see also~\cite{Duprez}). Moreover, this minimal time depends on the position of the open control set $\omega$ with respect to $\supp\, a_{12}$ (see~\cite{CherifMinimalTimeDisjointB}). In the second case, the null controllability result could depend on the coefficient $a_{11}$ and also on the position of $\omega$ (see~\cite{Duprezlissy1} and~\cite{Duprezlissy2}).}

Concerning the boundary controllability of systems of parabolic equations when $D$ is a positive multiple of the identity matrix, i.e., $D = d I_n$, with $d>0$, and $m < n$, the first results has been obtained in~\cite{fernandezcaraboundary2010} in the case $n = 2$ and $m=1$. In this paper, the authors prove that system~\eqref{eq:syst bound} is approximately and null controllable at time $T>0$ if and only if $A$ and $B$ satisfy the algebraic Kalman condition~\eqref{kalman} and 
	$$
\nu_1 - \nu_2 \not= d (\mu_j - \mu_ k) , \quad \forall k,j \in \N^\star \quad \hbox{ with } k \not= j,
	$$
where $\nu_1, \nu_2 \in \C^2$ are the eigenvalues of $A$ and $\mu_k := k^2$ are the eigenvalues of the Dirichlet-Laplace operator in $(0,\pi$). {The above condition shows the different nature of the null controllability problem for systems~\eqref{eq:syst int} and~\eqref{eq:syst bound}. When $D = d I_n$, with $d>0$, the null controllability of system~\eqref{eq:syst int} is equivalent to condition~\eqref{kalman spec} (which, in fact, is equivalent to~\eqref{kalman} when $D$ is a multiple of the identity matrix; see~\cite{AKBDGB09} or~\cite{AKBDGB09bis}).}

{The boundary null controllability of system~\eqref{eq:syst bound} when $D = d I_n$ ($d>0$)} has been generalized in~\cite{AKBGBdT11} to the case $n \ge 2$ and $m \in \N^\star$. The authors prove that system~\eqref{eq:syst bound} is null (resp., approximate) controllable at time $T$ if and only if
	\begin{equation*}
\hbox{rank}\, [\mathcal{L}_k  \, : \, \mathcal{B}_k]=nk~\hbox{ for all }k\in\mathbb{N^{\star}},
	\end{equation*}
where
	\begin{equation*}
\mathcal{B}_k:=\left[\begin{array}{c}B\\\vdots \\B\end{array}\right]
~\hbox{ and }
	\mathcal{L}_k=\left[
		\begin{array}{ccccc}
			 L_1	&	0	& 	\cdots	& 	0	\\
			 0	&	\ddots	& 	\ddots	&  	\vdots	\\
			 \vdots	&	\ddots	&	\ddots	& 	0	\\
			 0	&\cdots	& 	0  	& 	L_k
		\end{array}
	\right],
	\end{equation*}
with $L_k:=\mu_k I_n+A$ (see~\cite[Theorem $1.1$]{AKBGBdT11}). 

\begin{remark}\label{r1.1}
In particular, in~\cite{fernandezcaraboundary2010} and~\cite{AKBGBdT11}, the following property is proved ($D = d I_n$, with $d>0$): \textit{``Assume that condition~\eqref{kalman} holds (or equivalently, assume system~\eqref{eq:syst int} is null controllable at time $T>0$). Then, there exists a closed subspace $\mathcal{X}_0 \subset H^{-1}(0, \pi; \R^2)$, with finite codimension, which satisfies the following property: given $y^0 \in H^{-1}(0, \pi; \R^2)$, there exists a control $ v \in L^2 (0,T)$ such that the solution $y $ of system~\eqref{eq:syst bound} satisfies $y ( T, \cdot) = 0$ in $(0, \pi)$ if and only if $y^0 \in \mathcal{X}_0$.''}

Thus, we deduce that if system~\eqref{eq:syst int} is null controllable at time $T>0$, then system~\eqref{eq:syst bound} is also null controllable at time $T > 0$, apart from a finite-dimensional space.
\end{remark}

The case where $D$ is a diagonal matrix different from a multiple of the identity or the case in which the coupling matrices depend on the spatial variable are more delicate and new phenomena in the parabolic setting arise. For instance, a minimal time for the null controllability can appear, see \cite{AKBGBdT14}, and~\cite{CherifMinimalTimeDisjoint}. 


\paragraph{Non-diagonalizable diffusion matrices:}
In the case of non-diagonalizable diffusion matrices $D$, only partial results about controllability of systems~\eqref{eq:syst int} or~\eqref{eq:syst bound} have been established. In fact, {in the distributed control setting,} these results have been established in any spatial dimension for uniform elliptic time-independent operators $L$. More precisely:

\begin{itemize}
\item In~\cite{FCGBdT16} and under the following condition
	\begin{equation}\label{cond jordan}
\hbox{The dimensions of the \textit{Jordan blocks} of the canonical form of $D$  are }\leq 4,
	\end{equation}
the authors prove that system~\eqref{eq:syst int} is null (resp., approximately) controllable at time $T > 0$ with distributed controls if and only if~\eqref{kalman spec} holds. The technical condition \eqref{cond jordan} is a restriction due to the method used to provide the characterization (global Carleman estimates for scalar parabolic operators; see \cite[Remark $2.5$ and Section $5$]{FCGBdT16}). In particular, under condition~\eqref{cond jordan}, the authors prove a general result of approximate and null controllability at time $T>0$ for~\eqref{eq:syst int} and~\eqref{eq:syst bound} if $m \ge n$ and condition $\hbox{rank}\, B = n$ holds. This general result is also valid in the case of coupling matrices $A=A(t,x)$ which depend on $t$ and~$x$ or uniform parabolic operators $L = L(t)$ depending on $t$. 
\item In~\cite{lissy_zuazua}, the authors provide a complete answer for the problem of controllability of system~\eqref{eq:syst int} in the constant case without imposing any extra assumption on the Jordan blocks of $D$. In fact, they prove that system~\eqref{eq:syst int} is null (resp.,~approximately) controllable at time $T>0$ if and only if the constant matrices $D$, $A$ and $B$ satisfy the Kalman condition~\eqref{kalman spec} (see~\cite{lissy_zuazua} for more details). The approach followed in~\cite{lissy_zuazua} (the Lebeau-Robbiano strategy together with a precise study of the cost of controllability for linear ordinary differential equations) cannot be applied to system~\eqref{eq:syst bound}. However, as a consequence of their controllability results, it is not difficult to deduce that system ~\eqref{eq:syst bound} is approximately and null controllable at any time $T$ when $m \ge n$ and condition $\hbox{rank} \, B = n$ holds.
\item {Finally, in~\cite{GB-SN} the authors study the boundary null controllability of a one-dimensional phase field system of Caginalp type which is a model describing the transition between the solid and liquid phases in solidification/melting processes of a material occupying the  interval $(0, \pi)$. To this end, the authors prove the boundary null controllability of a linear $2 \times 2$ parabolic system with a non-diagonalizable diffusion matrix and a scalar control.}
\end{itemize}

The general goal of this paper is to study the controllability properties of coupled parabolic systems with non-diagonalizable diffusion matrices. To this end, we shall consider two examples of parabolic systems whose controllability properties cannot be obtained as a consequence of the controllability results proved in~\cite{FCGBdT16} and~\cite{lissy_zuazua}.

First, we will consider system~\eqref{eq:syst bound} when the matrices \MD{$D$, $A \in \mathcal{L}(\R^n)$} and $B\in \R^n$ ($n \ge 2$) are given by the expressions
\begin{equation}\label{def A B D}
	D=\left[
		\begin{array}{ccccc}
			 \displaystyle d	&	1	& 	0	&   	\cdots	& 	0	\\
			 0	&	d	& 	1	& 	\cdots  	& 	0	\\
			 \vdots	&	\vdots	&	\ddots	&  	\ddots	& 	\vdots	\\
			 0	&0 & 	\cdots	&  	d	& 	1	\\
			 0	&	0	& 	\cdots	& 	0  	& 	d
		\end{array}
	\right], \quad 
	A=\left[
		\begin{array}{cccc}
			 \displaystyle 0	&	0	& 		\cdots&   	0	\\
			 \vdots	&		\vdots	& \ddots & 	   	\vdots	\\
			0	& 0	& \cdots	&  	0		\\
			 \alpha &	0 & 	\cdots	&  	0	\
		\end{array}
	\right]\hbox{ and }
	B=\left[
		\begin{array}{c}
			 0		\\
			 \vdots		\\
			 0		\\
			 1
		\end{array}
	\right],
\end{equation}
for constants $d \ge 1$ and $\alpha \in \R$. Observe that since~\eqref{eq:syst bound} is a boundary controllability problem, its controllability properties cannot be deduced from~\cite{lissy_zuazua}. On the other hand, when $n \ge 5$, since $D$ does not satisfy condition~\eqref{cond jordan} we cannot apply the results in~\cite{FCGBdT16}, not even if $m \ge n$ and $\hbox{rank}\,B=n$.

\begin{remark}
It is not difficult to see that, when $d \ge 1$, the matrix $D$ is positive definite. We can conclude that systems~\eqref{eq:syst int} and~\eqref{eq:syst bound} are well-posed when, resp., $(y^0,u) \in L^2(0, \pi; \R^n) \times L^2(Q_T)$ and $(z^0,v) \in H^{-1}(0, \pi; \R^n) \times L^2(0,T)$.
\end{remark}

Secondly, we will consider problem~\eqref{eq:syst int} in the case in which $A$ is a matrix depending on $x$. To be precise, consider the following system
	\begin{equation}\label{syst simpl}
	\left\{\begin{array}{ll}
\partial_t y_1 -\partial_{xx} y_1 + q(x) y_1 = \partial_{xx} y_2 &\hbox{in } Q_T,\\
	\noalign{\smallskip}
\partial_t y_2 - \partial_{xx} y_2 ={1}_{\omega}u &\hbox{in } Q_T,\\
	\noalign{\smallskip}
y_1(\cdot,0)=y_2(\cdot,0)=y_1(\cdot,\pi)=y_2(\cdot,\pi)=0&\hbox{in } (0,T),\\
	\noalign{\smallskip}
y_1(0,\cdot)=y^0_1, \quad y_2(0,\cdot)=y^0_2&\hbox{in }(0,\pi),
	\end{array}\right.
	\end{equation}
where $q\in \mathcal{C}^{\infty}([0,\pi])$ is a given function, $y^0 = (y^0_1,y^0_2)^\star \in L^2(0,\pi;\mathbb{R}^2)$ is the initial datum and $u \in L^2(Q_T)$ is a scalar control.
	Observe that system~\eqref{syst simpl} has the same structure as system~\eqref{eq:syst int} with a coupling matrix $A$ depending on $x$. To be precise, $n = 2$, $m=1$, $D$ given in~\eqref{def A B D}, with $d = 1$, and $A = q(x) A_0$ with
	\begin{equation}\label{f5}
A_0 = \left[
	\begin{array}{cc}
1 & 0 \\ 0 & 0
	\end{array}
	\right] 
\quad \hbox{and} \quad
B=
	\left[
	\begin{array}{c}
			 0		\\
			 1
		\end{array}
	\right].
	\end{equation}
Again, $D$ is a definite positive matrix and, therefore, for any $y_1^0, y_2^0 \in L^2 (0, \pi)$ and $u \in L^2(Q_T)$, system~\eqref{syst simpl} has a unique solution $y \in L^2 (0,T; H_0^1(0, \pi; \R^2)) \cap \mathcal{C}^0([0,T]; L^2(0,\pi;\R^2))$ which depends continuously on the data.

\begin{remark}\label{r1.3}
When $n = 2$ and $q $ is a constant function, i.e., $q(x) = q_0$ for any $x \in (0, \pi)$, with $q_0 \in \R$, conditions~\eqref{kalman spec} and~\eqref{cond jordan} hold with $A = q_0 A_0$. In fact, the controllability matrix  $\left[ \mu_kD+A \, : \, B \right]$ does not depend on $q_0$. Thus, we can apply the results in~\cite{FCGBdT16} and~\cite{lissy_zuazua} and deduce that system~\eqref{syst simpl} is approximately and null controllable at any time $T>0$. Nevertheless, these results in~\cite{FCGBdT16} and~\cite{lissy_zuazua} cannot be applied when $A$ depends on $x$. To our knowledge, the controllability properties of system~\eqref{eq:syst int} with coupling matrices depending on $x$ are completely open. 
\end{remark}

Let us now present our first main result concerning system~\eqref{eq:syst bound}:
%


\begin{theorem}\label{thm:bord}
Let us consider the matrices $D, A \in \mathcal{L}(\R^n)$ and $ B \in \R^n$ given by~\eqref{def A B D}, with $d \ge 1$ and $\alpha \in \R$. When $\alpha >0$ we assume, in addition, that $n$ is odd. Then, 
\begin{enumerate}
\item If $\alpha \not= 0$, system~\eqref{eq:syst bound} is  null (resp., approximately) controllable at time $T$ if and only if the family of eigenvalues of the operator $L^\star:=- D^\star \partial_{xx} + A^\star $ has geometric multiplicity equal to $1$ and 
	\begin{equation}\label{f9}
\frac{\sqrt{\left| \alpha \right|}}{d^{ n/2}} \not\in \N^\star, 
	\end{equation}
if $n $ is odd and $ \alpha < 0$. 
\item When $\alpha =0$, system~\eqref{eq:syst bound} is approximately and null controllable at time $T > 0$.
\end{enumerate}
\end{theorem}

%
%


\begin{remark}\label{r1.4}
	Let us point out some consequences of Theorem~\ref{thm:bord}:
\begin{itemize}
\item Unlike system~\eqref{eq:syst bound}, system~\eqref{eq:syst int} is approximately and null controllable at time $T >0 $ (even in the $N$-dimensional case) when $D$, $A$ and $B$ are given by~\eqref{def A B D}. Indeed, it is easy to see that $\left[ \mu_kD+A \, : \, B \right]$ is a squared matrix and 
	$$
\left| \det \left[ \mu_kD+A \, : \, B \right] \right| = \mu_k \mu_k^2 \dots \mu_k^{n-1} \not= 0, \quad \forall k \ge 1,
	$$
($\mu_k = k^2 $, $k \ge 1$, are the eigenvalues of the Dirichlet-Laplace operator in $(0,\pi$)). Thus, condition~\eqref{kalman spec} holds and system~\eqref{eq:syst int} is approximately and null controllable at time~$T$. {Again, this shows the important differences between the null controllability properties of systems~\eqref{eq:syst int} and~\eqref{eq:syst bound}.}
\item We will see that, under the assumptions of Theorem~\ref{thm:bord}, the spectrum of the operator $L^\star :=- D^\star \partial_{xx} + A^\star $ is simple, apart from a finite number of eigenvalues (see Proposition~\ref{prop: index 0}). Therefore, Theorem~\ref{thm:bord} implies that system~\eqref{eq:syst bound} is always null controllable at any time $T>0$ apart from a finite dimensional space of $H^{-1}(0, \pi; \mathbb{C}^n)$ (in particular, this is the case when $n$ is odd).

\item	Condition~\eqref{f9} is related to the Fattorini-Hautus test for the operator $L^\star$ and is necessary in order to obtain the approximate controllability of system~\eqref{eq:syst bound} when $n$ is odd and the coefficient $\alpha $ of the matrix $A$  (see~\eqref{def A B D}) is negative. Therefore, this condition is also necessary for the null controllability of~\eqref{eq:syst bound}. No additional condition on $\alpha$ and $d$ is needed when $n$ is odd and $\alpha >0$. 

\item We will also see that, under assumptions of Theorem~\ref{thm:bord}, we can obtain a positive result on null controllability at any time $T>0$ for system
	\begin{equation}\label{f16}
	\left\{
	\begin{array}{ll}
y_t  -Dy_{xx} + Ay=B \delta_{x_0} u	&\hbox{ in } Q_T,\\
	\noalign{\smallskip}
y(\cdot,0)=y(\cdot,\pi)= 0 &\hbox{ in } (0,T),\\
	\noalign{\smallskip}
y(0,\cdot)= y^0 &\hbox{ in } (0,\pi),
	\end{array}
	\right.
	\end{equation}
where $u \in L^2(0,T)$, choosing very carefully the point $x_0 \in (0,\pi)$ where the control acts ($\delta_{x_0}$ is the Dirac distribution at point $x_0$). See Theorems~\ref{pointwise-app} and ~\ref{pointwise-null}.

\item	Observe that Theorem~\ref{thm:bord} does not provide any controllability result for system~\eqref{eq:syst bound} when $n = 2p$, $p \in \N^\star$, and $\alpha > 0$. In this case, the operator $L^\star :=- D^\star \partial_{xx}f + A^\star $ has, for any $k \ge 1$, exactly two real eigenvalues, given by
	$$
\lambda_{0,k} = dk^2 + \alpha^{1/n} k^{2 - 1/p} \quad \hbox{and} \quad \lambda_{p ,k} = dk^2 - \alpha^{1/n} k^{2 - 1/p } , 
	$$
and $n-2$ complex eigenvalues (see Proposition~\ref{eigenvectors}). In this case, the real eigenvalues of $L^\star$ could concentrate. As a consequence, the controllability problem for system~\eqref{eq:syst bound} could have a minimal time $T_0 \in [0,+ \infty]$ of null controllability which is related to the condensation index of the sequence $\left\{ \lambda_{0,k},\lambda_{p,k}\right\}_{k \ge 1}$ (see~\cite{AKBGBdT14} and Remark~\ref{r22}). In any case, condition~\eqref{f9} is also a necessary condition for the null controllability at time $T$ for system~\eqref{eq:syst bound} when $\alpha >0$ and $n$ is even {(see \eqref{eq:dalphaB*})}. The case $n=2$ is simpler and will be completely analyzed in Theorem~\ref{tn=2}.

\item The main advantage of the moment method, used in the present paper, is that it seems to be the best method to treat boundary null controllability problems with a reduced number of controls. Indeed, Theorem \ref{thm:bord} is the first result dealing with the boundary controllability of non-diagonalizable systems of parabolic equations when $n \ge 3$ (see also~\cite{GB-SN} where a similar problem is considered when $n=2$).
 \hfill $\Box$

\end{itemize}
\end{remark}

Theorem~\ref{thm:bord} provides sufficient conditions on $n$ and on the matrices $D, A \in \mathcal{L}(\R^n)$, given by~\eqref{def A B D}, which guarantee the approximate and null controllability of system~\eqref{eq:syst bound} at time $T>0$. As we already mentioned in the previous remark, this theorem does not cover the case $n = 2p$, with $p \in \N^\star$, and $\alpha >0$. In order to complete the study of the controllability problem for system~\eqref{eq:syst bound}, let us see the case $n = 2$ and $\alpha > 0$. One has:


\begin{theorem}\label{tn=2}
Assume that $ n = 2$ and let us consider the matrices $D, A \in \mathcal{L}(\R^2)$ and $ B \in \R^2$ given by~\eqref{def A B D}, with $d \ge 1$ and $\alpha > 0 $. Then, system~\eqref{eq:syst bound} is null (resp., approximately) controllable in $H^{-1}(0, \pi; \R^2)$ at time $T > 0$ if and only if
	$$
\frac{\sqrt \alpha}d \not\in \N^\star .
	$$
Moreover, if ${\sqrt \alpha}/d \in \N^\star$, there exists a closed subspace $\mathcal{X} \subset H^{-1}(0, \pi; \R^2)$, with infinite codimension, which satisfies the following property: given $y^0 \in H^{-1}(0, \pi; \R^2)$, there exists a control $ v \in L^2 (0,T)$ such that the solution $y $ of system~\eqref{eq:syst bound} satisfies $y ( T, \cdot) = 0$ in $(0, \pi)$ if and only if $y^0 \in \mathcal X$.
\end{theorem}


\begin{remark}\label{r1.5}
Even in the simplest case $n = 2$, Theorem~\ref{tn=2} shows an important difference with respect to the results on boundary controllability proved in~\cite{fernandezcaraboundary2010} and~\cite{AKBGBdT11}: When $D = d I_n$, if system~\eqref{eq:syst int} is null controllable at a time $T_0 > 0$, then system~\eqref{eq:syst bound} is also null controllable at any time $T > 0$, apart from a finite-dimensional space (see Remark~\ref{r1.1}). This property fails when $D$ is not diagonalizable ({even if $D$ is equal to a unique Jordan block}).
\end{remark}

In order to obtain Theorems~\ref{thm:bord} and~\ref{tn=2} we have used that the zero order coupling matrix $A$ in system~\eqref{eq:syst bound} is constant. Something similar occurs in~\cite{FCGBdT16} and~\cite{lissy_zuazua}: the authors use in a fundamental way that the zero order coupling matrix $A$ is constant. Let us now see that, if $A$ depends on $x$, the controllability properties of system~\eqref{eq:syst int} can be strongly affected. 	

Our third and last result is related to the controllability properties of system~\eqref{syst simpl}. It reads as follows:
\begin{theorem}\label{theo:negatif}

There exists a coefficient $q\in \mathcal{C}^{\infty}([0,\pi])$ such that:
\begin{enumerate}

\item [a)] There exists an open interval  $\omega\subset\subset (0,\pi)$ such that system \eqref{syst simpl} is never approximately controllable (then not null controllable)  for any time $T>0$.
\item[b)] There exists an open interval $\omega\subset\subset (0,\pi)$ such that system~\eqref{syst simpl} is null controllable (then approximately controllable) at any time $T>0$;
\end{enumerate}
\end{theorem}


%

\begin{remark}
{We can see system~\eqref{syst simpl} as a cascade system (see system~\eqref{cascade}) where the coupling term {$a_{12} y_2$} has been changed by the second order term {$a_{12}\partial_{xx} y_2$} and the coefficients are given by 
    $$
d_1 = d_2 = 1, \quad a_{11} = q, \quad a_{12} = - 1, \quad \text{and} \quad a_{21} = a_{22} = 0. 
    $$
Observe that this choice of coefficients implies condition~\eqref{cascade-cond} and the null controllability of system~\eqref{cascade} at any time $T >0$ and for any control open set $\omega \subset (0, \pi)$. Comparing the controllability result for system~\eqref{cascade} and system~\eqref{syst simpl}, Theorem~\ref{theo:negatif} shows that the controllability results for coupled parabolic systems with non-diagonalizable diffusion matrices may be very different and the location of the control domain plays a key role. A similar result to Theorem~\ref{theo:negatif} has been proved in~\cite{Duprezlissy2} when the coupling term in system~\eqref{syst simpl} is given by $\partial_x y_2$.
}
\end{remark}


\begin{remark}
{The proof of the Theorem~\ref{theo:negatif} relies on the construction of a non constant coefficient $q$ in $\mathcal{C}^{\infty}([0,\pi])$ which is constant in the interval $(5 \pi /12, 7 \pi/12)$. With this coefficient, we prove that system~\eqref{syst simpl} is not approximately controllable at any time $T>0$ when $\omega = (5 \pi /12, 7 \pi /12)$ and it is null controllable when we take $\omega \subset \supp q_x$ (i.e. $q$ is not constant in $\omega$). 
As stated in Remark~\ref{r1.3}, system~\eqref{syst simpl} is null controllable at any time $T >0$ when $q$ is a constant function.}
\end{remark}

\begin{remark}
{System~\eqref{syst simpl} is null controllable at time $T > 0$ for any open interval $\omega \subset (0, \pi)$ and $q \in L^\infty (0, \pi)$ if two independent distributed controls are exerted, one on the right-hand side of the first equation and one on the right-hand side of the second equation, see~\cite{FCGBdT16}. 
}
\end{remark}

%

\smallskip

{The rest of the paper is organized as follows. In Section \ref{sec:preliminar}, we study the spectrum and eigenvectors associated to systems \eqref{eq:syst int} and \eqref{eq:syst bound} and provide some properties needed to formulate the moment problem. In Section~\ref{sec:bound} we provide the proofs of Theorems~\ref{thm:bord} and \ref{tn=2}. Finally, in Section~\ref{sec_negative_} we prove Theorem \ref{theo:negatif}.}


\section{Preliminaries}\label{sec:preliminar}

\subsection{Spectral analysis}

{In the sequel, let us consider the following linear operator:
\[
	\begin{array}{rcl}
		L: D( L)\subset L^2(0,\pi;\mathbb{C}^n)&\longrightarrow& L^2(0,\pi;\mathbb{C}^n)\\\noalign{\smallskip}
 	 	f&\mapsto&Lf:=- D \partial_{xx}f + Af
	\end{array}
\]
	and its adjoint 
\[
	\begin{array}{rcl}
 		L^\star: D( L^\star)\subset L^2(0,\pi;\mathbb{C}^n)&\longrightarrow& L^2(0,\pi;\mathbb{C}^n)\\\noalign{\smallskip}
 		f&\mapsto&L^\star f:=- D^\star \partial_{xx}f + A^\star f, 
	\end{array}
\]
	where $D( L)=D(L^\star):= H^2(0,\pi;\mathbb{C}^n)\cap  H^1_0(0,\pi;\mathbb{C}^n)$ and  $D, A \in \mathcal{L}(\R^n)$ are given by~\eqref{def A B D}, with $n\geq2$, $d \ge 1$ and $\alpha \in \R$.
}


	It is well known that the operator $-\partial_{xx}:H^2(0,\pi;\mathbb{C})\cap  H^1_0(0,\pi;\mathbb{C})\longrightarrow L^2(0,\pi;\mathbb{C})$ admits a sequence of positive eigenvalues and a sequence of normalized eigenfunctions $\{w_k\}_{k\geq1}$, which is a Hilbert basis of $L^2(0,\pi)$, given by
	\begin{equation*}
\mu_k:=k^2\quad\hbox{and}\quad w_k(x):=\sqrt{\frac{2}{\pi}}\sin(kx), \quad x\in(0,\pi),~k\geq1.
	\end{equation*}

	Concerning the operator $L^\star$, we have the following description of its spectrum:
%
%
\begin{proposition}\label{eigenvectors}
Let us consider the matrices $D, A \in \mathcal{L}(\R^n)$ given by~\eqref{def A B D}, with $n \ge 2$, $d \ge 1$ and $\alpha \in \R$. Then, the following assertions hold:
{
	\begin{itemize} 
	\item [a)] If $\alpha\neq0$ then the spectrum of $L^\star$ is given by
$$
	\left\{
	\begin{array}{ll}
\sigma( L^\star)= \left\{\lambda_{j,k}:=d k^2 + \alpha^{1/n} k^{2-{\frac{2}{n}}}e^{{\frac{2\pi j}{ n}}i}:  k \in \N^\star, j \in \{0,...,n-1\} \right\},  &\hbox{if } \alpha > 0, \\
	\noalign{\smallskip}
\sigma( L^\star)= \left\{\lambda_{j,k}:=d k^2 + |\alpha|^{1/n} k^{2-{\frac{2}{n}}}e^{{\frac{(2 j + 1) \pi}{ n}}i}:  k\in\mathbb{N}^\star, j\in\{0,...,n-1\} \right\}, & \hbox{if } \alpha < 0.
	\end{array}
	\right.
	$$
	Moreover, for all $k\in \mathbb{N}^\star$ and $j\in \{0,\ldots,n-1\}$, an eigenvector of $L^\star$, associated to the eigenvalue $\lambda_{j,k}$, is given by
	\begin{equation*}
\Phi_{j,k} (x):=V_{j,k}w_k (x),
	\end{equation*}
where
	\begin{equation}\label{expr phi_k}
	V_{j,k}:=\frac{\left( c_{j,k}^l \right)_{1 \le l \le n}}{\left| \left( c_{j,k}^l \right)_{1 \le l \le n} \right|} \quad\hbox{with}\quad
c_{j,k}^l =
	\left\{
	\begin{array}{ll}
\left[\alpha^{-\frac 1n} k^{\frac{2}{ n}}e^{{\frac{-2\pi j}{ n}}i} \right]^{l-1}& \hbox{if } \alpha >0, \\
	\noalign{\smallskip}
\left[| \alpha |^{-\frac{1}n} k^{{\frac{2}{ n}}}e^{- \frac{(2 j + 1) \pi}n i}\right]^{l-1} & \hbox{if } \alpha <0.
	\end{array}
	\right.
	\end{equation}
\item [b)]If $\alpha=0$ then the spectrum of $L^\star$ is given by
$$
\sigma( L^\star)= \left\{\lambda_{k}:= 
d k^2 :  k\in\mathbb{N}^\star\right\}.
$$
Moreover, for all  $k\in\mathbb{N}^\star$,  $\Phi_{0,k} (x):=e_nw_k$ is an eigenvector of $L^\star$ associated to the eigenvalue $\lambda_{k}$ and $\Phi_{1,k} (x):=e_{n-1}w_k$, \dots, $\Phi_{n-1,k} (x):=e_{1}w_k$ are generalized eigenvectors of $L^\star$ associated to $\lambda_{k}$.
	\end{itemize}
}
\end{proposition}

\begin{proof}
\textbf{Case} $\boldsymbol{\alpha\neq0}$:
	The goal is to solve the following eigenvalue problem:
\[
	L^\star\Psi=\lambda \Psi, \quad \lambda \in \C, \quad {\Psi \in H^2(0,\pi;\mathbb{C}^n)\cap  H^1_0(0,\pi;\mathbb{C}^n)}.
\]
	Considering $\Psi=(\psi_1,\ldots,\psi_n)$, the previous eigenvalue problem is equivalent to
	\[
	\left\{
	\begin{array}{l}
-d\partial_{xx}\psi_1+ \alpha \psi_n=\lambda\psi_1\\
	\noalign{\smallskip}
-\partial_{xx}\psi_1-d\partial_{xx}\psi_2 =\lambda\psi_2\\
	\noalign{\smallskip}
\qquad \qquad \vdots  \\
	\noalign{\smallskip}
-\partial_{xx}\psi_{n-1}-d\partial_{xx}\psi_n =\lambda\psi_n.
	\end{array}
	\right.
	\]
	It can be rewritten as an algebraic eigenvalue problem 
	\[
	\left\{
	\begin{array}{l}
dk^2c_{k}^1+ \alpha c_{k}^n = \lambda c_{k}^1\\
	\noalign{\smallskip}
k^2c_{k}^1 + dk^2c_k^2 = \lambda c_{k}^2\\
	\noalign{\smallskip}
\qquad \qquad \vdots  \\
	\noalign{\smallskip}
k^2c_k^{n-1}+dk^2c_k^{n}=\lambda c_k^{n},
	\end{array}
	\right.
	\]
for all $k\geq1$. From the previous expression, we obtain
	$$
	\left\{
	\begin{array}{l}
(\lambda - dk^2 ) c_{k}^1 = \alpha c_{k}^n, \\
	\noalign{\smallskip}
\left(\lambda - dk^2 \right) c_{k}^l = k^2 c_{k}^{l-1}, \quad \forall l = 2, \dots , n.
	\end{array}
	\right.
	$$
After some computations, we also get
	\begin{equation}\label{expr phi_kbis}
c_{k}^l = k^{2(l -1)} \left(\lambda - dk^2 \right)^{-(l-1)} c_{k}^1, \quad \forall l = 1, 2, \dots , n.
	\end{equation}
	Due to the identity for $l=n$, we necessarily have to impose the following condition
\[
	\alpha k^{2(n-1)}=(\lambda-dk^2)^{n},
\]
determining all the eigenvalues for the eigenvalue problem at the beginning of the proof. Taking into account that $\alpha \not=0$, we deduce that the previous equation has $n$ distinct solutions $\lambda_{j,k}$, $j= 0, \dots, n-1$, which are given explicitly in item $a)$. 

Finally, for all $k\geq 1$ and $j=0,\ldots,n-1$, from~\eqref{expr phi_kbis} and the expression of $\lambda_{j,k}$, it is not difficult to see that $\Phi_{j,k} $, with $c_{j,k}^l$ given in~\eqref{expr phi_k}, is an eigenvector of $L^\star$ associated to the eigenvalue $\lambda_{j,k}$. 

\smallskip

\textbf{Case} $\boldsymbol{\alpha=0}$:
The eigenvalue problem can be rewritten as the algebraic eigenvalue problem
	\[
	\left\{
	\begin{array}{l}
dk^2c_{k}^1= \lambda c_{k}^1\\
	\noalign{\smallskip}
k^2c_{k}^1 + dk^2c_k^2 = \lambda c_{k}^2\\
	\noalign{\smallskip}
\qquad \qquad \vdots \\
	\noalign{\smallskip}
k^2c_k^{n-1}+dk^2c_k^{n}=\lambda c_k^{n} \, .
	\end{array}
	\right.
	\]
We deduce that the previous problem has a unique eigenvalue $\lambda_k:=dk^2$ with algebraic multiplicity equal to $n$.  An associated eigenvector is the vector $e_n$. Moreover the associated generalized eigenvectors are $e_1, \dots , e_{n-1}$.
This ends the proof.
\end{proof}


	Let us consider the set 
	\begin{equation}\label{basis}
\mathcal{B}^\star := \left\{\Phi_{j,k} :\,k\in \mathbb{N}^\star,j\in\{0,...,n-1\} \right\},
	\end{equation}
where the functions $\Phi_{j,k} $ are given in Proposition~\ref{eigenvectors}.
	Then, we obtain the following result:
	%
	%
\begin{lemma}\label{riesz_bases}
{Under the assumptions of Proposition~\ref{eigenvectors}, the set $\mathcal{B}^\star$ is a Schauder basis of the spaces $L^2(0,\pi;\mathbb{C}^n)$ and $H^1_0(0,\pi;\mathbb{C}^n)$, normalized in $L^2(0,\pi;\mathbb{C}^n)$.}
\end{lemma}
\begin{proof}

	Let us prove the result in the case $\alpha >0$. The case $ \alpha < 0$ can be deduced with a similar reasoning and the case $\alpha=0$ is trivial.
	Consider the Schauder basis $\mathcal{B}_c$ of $L^2(0,\pi;\mathbb{C}^n)$ given by
\[
	\mathcal{B}_{c}:= \left\{e_iw_k:\,k\in \mathbb{N}^\star, \ i\in\{1,...,n\} \right\}.
\]   
	Let $h\in L^2(0,\pi;\mathbb{C}^n)$ (or $h\in H^1_0(0,\pi;\mathbb{C}^n)$). There exists a unique real sequence 
	$\{\alpha_{i,k}\}_{k\in\mathbb{N}^\star,1\leq i\leq n}$ such that
\begin{equation*}
			h=\sum\limits_{k=1}\limits^{+\infty}\sum\limits_{i=1}\limits^{n} \alpha_{i,k} e_i w_k \quad \hbox{in } L^2(0,\pi;\mathbb{C}^n) \quad \hbox{(resp., in } H^1_0(0,\pi;\mathbb{C}^n) ).
\end{equation*}
	{We remark that the matrix $\mathcal{V}_k:=(\widetilde V_{0,k}|\cdots|\widetilde V_{n-1,k})$ is a Vandermonde matrix (and so, it is invertible), where $\widetilde V_{j,k}= \left( c_{j,k}^l \right)_{1\leq l\leq n}$ for any $k\in \mathbb{N}^\star$ and $j\in\{0,\ldots,n-1\}$, see \eqref{expr phi_k}. }
	Therefore, for each $k\in\mathbb{N}^\star$, there exist unique $\gamma_{0,k},\ldots,\gamma_{n-1,k}\in\mathbb{C}^\star$ 
	such that 
\[
	\sum\limits_{i=1}\limits^{n} \alpha_{i,k} e_i=
	\sum\limits_{j=0}\limits^{n-1} \gamma_{j,k} V_{j,k}.
\]

	Finally, arguing by contradiction, we can obtain the sequence $\{\gamma_{p,k}\}_{k\in\mathbb{N}^\star,0\leq p\leq n-1}$ is unique.
\end{proof}

\begin{remark}
    In the sequel, we will use the notation 
    $$
        \Lambda_0:= \left\{\lambda_{j,k}: k \in \N^\star, j \in \{0,...,n-1\}\right\},
    $$
    where $\lambda_{j,k}$ is given in Proposition~\ref{eigenvectors}.
\end{remark}

\subsection{Biorthogonal family}

Given a complex sequence $\Lambda = \{\Lambda_k \}_{k\geq 1} \subset \mathbb{C}$, let us denote by $p_{k}$ the complex function given by 
	\begin{equation}\label{f0}
p_{k}(t):=e^{-\Lambda_{k}t}, \quad \forall t\in(0,T). 
	\end{equation}
{We will see in Section~\ref{sec:bound} that the existence of a biorthogonal family in $L^2 (0,T; \C)$ to the sequence $\{p_{k}\}_{k \ge 1}$ will play a key role in the study of the controllability of systems \eqref{eq:syst bound} and \eqref{f16}. }
Recall that the sequence $\{q_{k}\}_{k\geq1}$ is a biorthogonal family to the sequence $\{p_k \}_{k\geq1}$ in $L^2 (0,T; \C)$ if 
	\begin{equation*}
\displaystyle\int_0^T e^{-\Lambda_k t} q^\star_j(t) \, dt = \delta_{kj},\quad \forall k,j\geq1.
	\end{equation*}
One has:

\begin{theorem}\label{lemma: base ortho}
	Let $T>0$ and consider a sequence $\Lambda = \{\Lambda_k \}_{k\geq 1} \subset \mathbb{C}$
 satisfying
	\begin{equation}\label{H:lambda}
	\left\{\begin{array}{l}
\displaystyle \Lambda_i\neq\Lambda_k , \quad \forall i,k\in\mathbb{N}^{\star}\hbox{ with }i\neq k,\\
	\noalign{\smallskip}
\displaystyle \sum_{k\geq 1}\dfrac{1}{|\Lambda_k|}<+\infty  \quad \hbox{and} \quad \Re( \Lambda_k ) \geq \delta |\Lambda_k |> 0, \quad \forall k\geq 1,
	\end{array}
	\right.
	\end{equation}
for some positive constant $\delta$. 
	Then, there exists a biorthogonal family $\{q_{k}\}_{k\geq1}$ in $L^2 (0,T; \C)$ to the family $\{p_k \}_{k\geq1}$, given in~\eqref{f0}.
	Moreover, for any $\varepsilon>0$, there exists a constant $C_\varepsilon > 0$ 
	such that
	\begin{equation}\label{estim q_k}
\|q_{k}\|_{L^2(0,T; \C)}\leq C_\varepsilon e^{(c(\Lambda)+\varepsilon)\Re(\Lambda_{k})},
	\end{equation}
where $c(\Lambda) \in [0, +\infty]$ is the condensation index of the sequence $\Lambda$. 
\end{theorem}
	This result corresponds to \cite[Proposition $4.1$ and Remark 4.3]{AKBGBdT14}. See
	\cite[Definition $3.1$]{AKBGBdT14}, for the definition of the condensation index.


\


In some situations, this index of condensation of the sequence $\Lambda = \{\Lambda_k \}_{k\geq 1} \subset \mathbb{C}$ can be equal to zero. Let us consider two different situations:
\begin{lemma}\label{l2}
Let us consider a sequence $\Lambda = \{\Lambda_k \}_{k\geq 1} \subset \mathbb{C}$ satisfying \eqref{H:lambda}. Let us assume that there exist a positive constant $\rho > 0$ and a positive integer $k_1$ such that one of the following conditions hold

	\begin{equation}\label{cond lambda0}
\left| \Lambda_k - \Lambda_l \right| \geq \rho  |\Lambda_k |^{ 1/2} ,
\quad \forall k \geq k_1\hbox{  and }l \neq k.
	\end{equation}
	\begin{equation}\label{H:gap}
\left| \Lambda_{k}-\Lambda_{l} \right| \ge \rho |k - l|, \quad \forall k,l \in\mathbb{N}^\star . 
	\end{equation}
Then, 
	$$
c(\Lambda) = 0.
	$$
\end{lemma}
	For a proof of the previous lemma, we refer to \cite[Theorem $6$]{SHACKELL}.

	\smallskip

{In Section~\ref{sec:bound} we will also use a result on \MD{the} existence of biorthogonal families to some complex matrix exponentials. In order to state the result, let us fix $\eta \ge 1$, a positive integer, and let us introduce the notation: 
	\begin{equation*}
p_{k}^{(j)}(t) := t^j e^{- \Lambda_k t}, \quad \forall t >0, \quad (k \ge 1 \mbox{ and } j : 0 \le j \le \eta -1),
	\end{equation*}
where $\Lambda =\left\{ \Lambda _{k}\right\}_{k\geq 1}\subset \mathbb{C}$ is a sequence of complex numbers. 

Let us recall that the family $\{q_{k}^{(j)}\}_{k \ge 1, 0\le j\le \eta-1} \subset L^2(0, T; \C)$ is biorthogonal to the sequence $\{p_{k}^{(j)}\}_{k \ge 1, 0 \le j\le \eta-1}$ if the equalities 
	\begin{equation} \label{ortho}
\displaystyle \int_{0}^{T } t^{j}e^{-\Lambda _{k}t} q_{m}^{{(l)}\star}(t) \, dt = \delta_{km} \delta_{jl}, \quad \forall (j,k), (l, m) : k, m\geq 1, \hbox{ } 0 \leq  j, l \leq \eta - 1,
	\end{equation}
hold.

With the previous notation, one has:
\begin{theorem}\label{t3} 
Let us fix $\eta \ge 1$, a positive integer, and $T >0$. Assume that $\left\{ \Lambda _{k}\right\}_{k\geq 1}$ is a sequence
of complex numbers satisfying~\eqref{H:lambda} and~\eqref{H:gap} for two positive constants $\delta$ and $\rho$. Then, there exists a family $\left\{ q_{k}^{(j)} \right\}_{k\geq 1, 0 \leq j \leq \eta - 1}\subset L^2 (0,T; \C)$ biorthogonal to $\left\{p_{k}^{(j)} \right\}_{k\geq 1, 0 \leq j \leq \eta - 1}$ such that, for every $\varepsilon >0$, there exists $ C_\varepsilon > 0$ for which 
\begin{equation}  \label{biacot}
	\| q_{k}^{(j)} \|_{L^2(0, T ; \C )} \le C_\varepsilon e^{\varepsilon \Re \left( \Lambda_k \right)}, \quad \forall (j, k) : k \ge 1, \ 0 \le j \le \eta - 1.
\end{equation}
\end{theorem}

For a proof of this result, see~\cite[Theorem 1.2]{AKBGBdT11}.
}
	
\subsection{Condensation index of the sequence $\Lambda_0$}	
Our next objective is to check conditions~\eqref{H:lambda} and \eqref{cond lambda0} for the sequence 
    $$
\Lambda_0=\{\lambda_{j,k}\}_{0 \le j \le n-1,k \ge 1}
    $$ 
in the case $\alpha \not= 0$. One has:

\begin{proposition}\label{prop: index 0}
Let us consider the matrices $D, A \in \mathcal{L}(\R^n)$ given by~\eqref{def A B D}, with $n \ge 2$, $d \ge 1$ and $\alpha \in \R$, $\alpha \not=0$. In addition, assume that $n:=2p+1$ ($p\in\mathbb{N}^{\star}$) when $\alpha >0$. Then, there exist a positive integer $k_1$ and a constant $\rho >0$, only depending on $n$, $\alpha$ and $d$, such that 
	\begin{equation}\label{cond lambda}
\left| \lambda_{j,k}-\lambda_{j',k'} \right| \geq \rho \left| \lambda_{j,k} \right|^{1/2}, \ \forall k \ge k_1, \ \forall k' \ge 1, \ \forall j,j' \in\{0,...,n-1\} :  (j,k)  \not= (j',k') .
	\end{equation}
The expression of $\lambda_{j,k} $ ($ k \ge 1$, $j: 0 \le j \le n-1$) is given in item a) of Proposition~\ref{eigenvectors}.
\end{proposition}

\begin{proof}
Let us prove that condition~\eqref{cond lambda} holds when $\alpha \not= 0$. Fix $(j,k) \in \{0,...,n-1\} \times \mathbb{N}^\star$ and notice that
	$$
\left| \lambda_{j,k} \right|^{1/2} \leq k\left( d + | \alpha |^{\frac 1n} k^{-\frac{2}{n}} \right)^{1/2} \leq  C_1k,
	$$
for all $k\geq 1$ and for some positive constant $C_1$ independent of $j$ and $k$. 

The goal is to prove that there exist a positive integer $k_1$ and a constant $C_2 > 0$, only depending on $n$, $\alpha$ and $d$, such that 
	\begin{equation}\label{distint}
\left| \lambda_{j,k}-\lambda_{j',k'} \right| \geq C_2 k, \quad \forall k \ge k_1, \quad \forall k' \ge 1, \quad \forall j,j' \in\{0,...,n-1\} :  (j,k)  \not= (j',k') .
	\end{equation}
Observe that condition~\eqref{cond lambda} is a direct consequence of two previous inequality~\eqref{distint}.

\medskip

\noindent\textbf{Case A: $\alpha >0$}. In this case, we assume that $n = 2p+1$, with $p \ge 1$. 

\medskip

Thanks to the fact that the function $ dr^2 + \alpha^{\frac 1n}\cos(\frac{2\pi j}{n})r^{2-\frac{2}{n}} $ is increasing for $r$ large enough, we deduce there exists $k_0 \in \mathbb{N}$ such that for any $k\geq k_0$, $k'\in \mathbb{N}^\star$ (with $k'\neq k$) and $j\in \{0,\ldots,n-1\}$, we have (see item a) of Proposition~\ref{eigenvectors})
	$$
	\begin{array}{l}
\left|\lambda_{j,k} - \lambda_{j,k'} \right| \geq \min \left\{ \left| \Re (\lambda_{j,k}) - \Re(\lambda_{j,k-1}) \right|, \left|\Re(\lambda_{j,k}) - \Re(\lambda_{j,k+1}) \right| \right\},\\ 
	\noalign{\smallskip}
\phantom{ \left|\lambda_{j,k} - \lambda_{j,k'} \right| } \geq C_3 k,
	\end{array}
	$$
where $C_3 > 0$ is independent of $j$, $k$ and $k'$. 
On the other hand, using that $n$ is odd, we can also prove that if $k \not= k'$ we have $\lambda_{j,k} \not= \lambda_{j,k'} $. This proves conditions~\eqref{distint} and~\eqref{cond lambda} for $k_1 =1$.

	Consider now $k\geq 1$ and $j,j'\in\{0,...,n-1\}$ such that $j'\neq j$.
	The goal is to prove inequality~\eqref{distint} for any $ k \ge k_1$, $k' \ge 1$ and for any $ j,j' \in\{0,...,n-1\}$ with $ j \not=j' $. Again, condition~\eqref{cond lambda} will be a direct consequence of inequality~\eqref{distint}.

	Let us introduce the notation 
	$$
c_j := \cos \left( 2\pi j/n \right) \hbox{ and } s_{j}:=\sin \left(2\pi j/n \right), \quad \forall j \in\{0,...,n-1\}.
	$$
Observe that thanks to the assumption $n =2p+1$, with $p \in \mathbb{N}^\star$, we can conclude that $s_j \not= 0$ for any $j : 1 \le j \le n-1$. On the other hand, with the previous notation, one has (recall that $ \alpha >0$)
	$$
\Re \left( \lambda_{j,k} \right) = dk^2 + \alpha^{\frac 1n} c_j k^{2 - \frac 2n} \hbox{ and } \Im \left( \lambda_{j,k} \right) = \alpha^{\frac 1n} s_j k^{2 - \frac 2n} , \quad \forall k \ge 1, \quad \forall j \in\{0,...,n-1\} .
	$$

In order to show~\eqref{distint}, we distinguish three cases:

\smallskip

\noindent\textbf{Case A.1:} $j=0$. 
 
\begin{enumerate}
\item [$i.$] If $k'>k/2$ and $j' \not= 0$, we obtain 
	$$
\left| \lambda_{0,k}-\lambda_{j',k'} \right| \geq \left| \Im \left( \lambda_{0,k}-\lambda_{j',k'} \right) \right| = \alpha^{\frac 1n} \left| s_{j'} \right| (k')^{2-{\frac{2}{n}}} \geq \alpha^{\frac 1n} \left| s_{j'} \right| k'\geq \frac{ 1}{2}\alpha^{\frac 1n} \left| s_{j'} \right| k;
	$$
\item [$ii.$] If $k'\leq k/2$ then we deduce 
	$$
\left| \lambda_{0,k}-\lambda_{j',k'} \right| \geq \left| d(k^2-(k')^2)+\alpha^{\frac 1n} k^{2-{\frac{2}{n}}}- \alpha^{\frac 1n} c_{j'}(k')^{2-{\frac{2}{n}}} \right| \geq d(k^2-(k')^2)\geq  \frac34 dk. 
	$$
In both cases,~\eqref{distint} holds for $k_1 = 1$.
\end{enumerate}

\smallskip

\noindent \textbf{Case A.2:} $j\neq0$ and $j' \not= j$ with $\text{sign}(s_j)\neq \text{sign}(s_{j'})$. Then, we have 
	$$
\left| \lambda_{j,k}-\lambda_{j',k'} \right| \geq  \left| \Im \left( \lambda_{j,k}-\lambda_{j',k'} \right) \right| \ge \alpha^{\frac 1n} \left| s_j \right| k^{2-{\frac{2}{n}}}\geq \alpha^{\frac 1n} \left| s_j \right| k;
	$$
In particular, we deduce~\eqref{distint} for $k_1 = 1$.

\smallskip

\noindent\textbf{Case A.3:} 
	$j\neq0$ and $j' \not= j$ with $\text{sign}(s_j)= \text{sign}(s_{j'})$.
	First, notice that $j'\neq j$ and the fact that $n$ is odd implies $s_j \not=0$, $s_{j'} \not=0$ and $|s_j|\neq|s_{j'}|$. Now, we define 
	$$
\beta := \left\{
	\begin{array}{ll}
\displaystyle \left(\frac{ \left| s_j \right| }{\left| s_{j'} \right| }\right)^{1/2}&\quad\mbox{ if }\quad\dfrac{\left| s_j \right|}{\left| s_{j'} \right|} > 1,\\
	\noalign{\smallskip}
\displaystyle \left(\frac{\left| s_{j'} \right|}{\left| s_j \right|}\right)^{1/2}&\quad\mbox{ if } \quad \dfrac{\left| s_j \right|}{\left| s_{j'} \right|} < 1
	\end{array}
	\right.
	$$
and the numbers:
	$$
\beta_1 := \left( \frac 1 \beta \dfrac{\left| s_j \right|}{\left| s_{j'} \right|} \right)^{\frac n{2(n-1)}}\quad\mbox{ and }\quad
\beta_2 := \left( \beta \dfrac{\left| s_j \right|}{\left| s_{j'} \right|} \right)^{\frac n{2(n-1)}} .
	$$
%
In any case, one has $\beta > 1$ and $\beta_1 < \beta_2$. 

Then, we consider three different cases:
\begin{enumerate}
\item [$i.$] If $k'\leq \beta_1k$, with $k \ge 1$, we use the imaginary part of $ \lambda_{j,k} - \lambda_{j',k'}$ to deduce 
	$$
	\begin{array}{l}
\displaystyle \left| \lambda_{j,k} - \lambda_{j',k'} \right| \ge \alpha^{\frac 1n} \left( k^{2-{\frac{2}{n}}} |s_j| - (k')^{2-{\frac{2}{n}}} |s_{j'}| \right) \ge \alpha^{\frac 1n} k^{2-{\frac{2}{n}}} \left(  |s_j| -  \beta_1^{\frac{2(n-1)}{n}} |s_{j'}| \right) \\
	\noalign{\smallskip}
\displaystyle \phantom{\left| \lambda_{j,k} - \lambda_{j',k'} \right|}  = \alpha^{\frac 1n}k^{2-{\frac{2}{n}}} \left(1-\frac{1}{\beta}\right) \left|s_j \right|  \geq \alpha^{\frac 1n} \left(1-\frac{1}{\beta}\right) \left|s_j \right| k.
	\end{array}
	$$
\item [$ii.$] Let us now assume that $k'\geq \beta_2k$.  Then a similar reasoning as before provides 
	$$
\left| \lambda_{j,k} - \lambda_{j',k'} \right| \geq  \alpha^{\frac 1n} (\beta-1) \left| s_j \right| k^{2-{\frac{2}{n}}} \geq  \alpha^{\frac 1n} (\beta-1) \left| s_j \right|k.
	$$
In particular we have~\eqref{distint} for $k_1 =1$.
\item [$iii.$] If now $k'=\gamma k$ with $\gamma\in(\beta_1,\beta_2)$, then we deduce
	$$
	\begin{array}{l}
\left| \lambda_{j,k}-\lambda_{j',k'} \right| \geq \left| d\,k^2+ \alpha^{\frac 1n} c_jk^{2-{\frac{2}{n}}}-d\,\gamma^2k^2- \alpha^{\frac 1n} c_{j'}\gamma^{2-{\frac{2}{n}}}k^{2-{\frac{2}{n}}} \right|\\
	\noalign{\smallskip}
\phantom{\left| \lambda_{j,k}-\lambda_{j',k'} \right|} = k^{2-{\frac{2}{n}}} \left| d(1-\,\gamma^2)k^{2/n}+ \alpha^{\frac 1n} (c_j-c_{j'}\gamma^{2-{\frac{2}{n}}}) \right|.
	\end{array}
	$$
	\end{enumerate}

Since 
	$$
	\left\{
	\begin{array}{ll}
\displaystyle \gamma^{2-2/n}\in \left(\frac{ \left|s_j \right|^{1/2}}{ \left|s_{j'} \right|^{1/2}},\frac{ \left| s_j \right|^{3/2}}{\left|s_{j'} \right|^{3/2}}\right)
&\mbox{ if }\dfrac{ \left|s_j \right| }{\left| s_{j'} \right|} > 1,\\
	\noalign{\smallskip}
\displaystyle \gamma^{2-2/n} \in \left( \frac{ \left| s_j \right|^{3/2}}{\left|s_{j'} \right|^{3/2}} ,\frac{ \left|s_j \right|^{1/2}}{ \left|s_{j'} \right|^{1/2}} \right)
&\mbox{ if } \dfrac{ \left|s_j \right| }{\left| s_{j'} \right|} < 1,
\end{array}\right.
	$$
we easily deduce that $\left| \gamma^2-1 \right| \geq C>0$ and then there exists a new positive integer $k_1\in \mathbb{N}^\star$ such that for any $k\geq k_1$ we have
	$$
\left| d(1-\gamma^2)k^{2/n}+ \alpha^{\frac 1n} (c_j-c_{j'}\gamma^{2-{\frac{2}{n}}}) \right| \geq C > 0.
	$$ 
Therefore, for $k\geq k_1$, we obtain
	$$
\left| \lambda_{j,k}-\lambda_{j',k'} \right| \geq Ck^{2-{\frac{2}{n}}}\geq Ck.
	$$

In conclusion, we have proved inequality~\eqref{distint} when $n$ is odd and $\alpha > 0$. 

\ 

\medskip

\textbf{Case B: $\alpha < 0$}. 

\medskip

In this case, the eigenvalues of the operator $L^\star = - D^\star \partial_{xx} + A^\star$ are given by
	$$
\lambda_{j,k} = d\,k^2 + |\alpha|^{1/n} k^{2-{\frac{2}{n}}}e^{{\frac{(2 j + 1) \pi}{ n}}i}, \quad  k\ge 1, \quad j\in\{0,...,n-1\}.
	$$

Observe that, in this case,  $\lambda_{p ,k} \in \R$ when $n = 2p + 1$ ($p \ge 1$, an integer) and $\lambda_{j,k} \in \C \setminus \R$, otherwise. 
Again, our goal is to prove inequality~\eqref{distint} for a positive integer $k_1$ and a positive constant $C_2$ only depending on $n$, $d$ and $\alpha$. 

	Let us introduce the notation 
	$$
\widetilde c_j := \cos \left( \frac{(2 j + 1) \pi}{n} \right) \quad \hbox{and} \quad \widetilde s_{j}:=\sin \left( \frac{(2 j + 1) \pi}{n} \right) , \quad \forall j \in\{0,...,n-1\}.
	$$
With this notation, one has (recall that $ \alpha < 0$)
	$$
\Re \left( \lambda_{j,k} \right) = dk^2 + | \alpha |^{\frac 1n} \widetilde c_j k^{2 - \frac 2n} \hbox{ and } \Im \left( \lambda_{j,k} \right) = |\alpha |^{\frac 1n} \widetilde s_j k^{2 - \frac 2n} , \quad \forall k \ge 1, \quad \forall j \in\{0,...,n-1\} .
	$$

Using the fact that the function $ dr^2 + | \alpha |^{\frac 1n}\widetilde c_j r^{2-\frac{2}{n}} $ is increasing for $r$ large enough, we deduce the existence of~$\widetilde k_0 \in \mathbb{N}$ such that for any $k\geq \widetilde k_0$, $k'\in \mathbb{N}^\star$ (with $k'\neq k$) and $j\in \{0,\ldots,n-1\}$, we have 
	$$
	\begin{array}{l}
\left|\lambda_{j,k} - \lambda_{j,k'} \right| \geq \min \left\{ \left| \Re (\lambda_{j,k}) - \Re(\lambda_{j,k-1}) \right|, \left|\Re(\lambda_{j,k}) - \Re(\lambda_{j,k+1}) \right| \right\},\\ 
	\noalign{\smallskip}
\phantom{ \left|\lambda_{j,k} - \lambda_{j,k'} \right| } \geq \widetilde C_3 k,
	\end{array}
	$$
where $\widetilde C_3 > 0$ is independent of $j$, $k$ and $k'$. 

As in the case $\alpha > 0$, let us consider $k\geq 1$ and $j,j'\in\{0,...,n-1\}$ such that $j'\neq j$ and let us prove inequality~\eqref{distint}. 
We distinguish four cases:

\smallskip

\noindent\textbf{Case B.1:} $\lambda_{j,k} \in \R$, i.e., $n = 2 p +1$ ($p \in \N^\star$) and $j=p$. In this case, $\lambda_{p, k} = dk^2 - | \alpha |^{\frac 1n} k^{2 - \frac 2n}$. So,
 
\begin{enumerate}
\item [$i.$] If $k'>k/2$ and $j' \not= p$, we obtain 
	$$
\left| \lambda_{p,k}-\lambda_{j',k'} \right| \geq \left| \Im \left( \lambda_{p,k}-\lambda_{j',k'} \right) \right| = |\alpha |^{\frac 1n} \left| \widetilde s_{j'} \right| (k')^{2-{\frac{2}{n}}} \geq |\alpha |^{\frac 1n} \left| \widetilde s_{j'} \right| k'\geq \frac{ 1}{2}|\alpha |^{\frac 1n} \left| \widetilde s_{j'} \right| k;
	$$
\item [$ii.$] If $k'\leq k/2$, we deduce 
	$$
\left| \lambda_{p,k}-\lambda_{j',k'} \right| \geq d[k^2-(k')^2] - | \alpha |^{\frac 1n} k^{2-{\frac{2}{n}}}- | \alpha |^{\frac 1n} \left| \widetilde c_{j'} \right| (k')^{2-{\frac{2}{n}}} \geq k^2 \left( \frac34 d  - 2 | \alpha |^{\frac 1n} k^{- \frac{2}{n}} \right).
	$$
From this expression, we deduce~\eqref{distint} for $k_1 \ge 1$ and $C_2 >0$ only depending on $n$, $d$ and~$\alpha$. 
\end{enumerate}

\smallskip

\noindent\textbf{Case B.2:} $\lambda_{j,k} \in \C \setminus \R$ and $\lambda_{j' ,k' } $ is such that $\text{sign} (\widetilde s_j) \not = \text{sign} (\widetilde s_{j'})$. In this case, $\widetilde s_j \not= 0$ and we can write
	$$
\left| \lambda_{j,k}-\lambda_{j',k'} \right| \geq  \left| \Im \left( \lambda_{j,k}-\lambda_{j',k'} \right) \right| \ge \left| \alpha \right|^{\frac 1n} \left| \widetilde s_j \right| k^{2-{\frac{2}{n}}}\geq \left| \alpha \right|^{\frac 1n} \left| \widetilde s_j \right| k;
	$$
In particular, we deduce~\eqref{distint} for $k_1 = 1$.

\smallskip

\noindent\textbf{Case B.3:} $\lambda_{j,k} \in \C \setminus \R$ and $\lambda_{j' ,k' } $ is such that $\widetilde s_j = \widetilde s_{j'}$, with $j \not= j'$. In this case, notice that  $\widetilde c_j = - \widetilde c_{j'} \not= 0 $, $\widetilde s_j \not= 0$ and $n$ should be an even number greater or equal than $4$. On the other hand, we also have:
	\begin{equation}\label{xM}
(x + M )^{2 - \frac 2n} - x^{2 - \frac 2n} \ge \left(2 - \frac 2n \right) M x^{1 - \frac 2n}, \quad \forall x,M>0, \quad (n \ge 2).
	\end{equation}
Let us fix $\eta$, the positive root of the equation $d x^2 + 2 d x - | \alpha |^{\frac 1n} \left| \widetilde c_{j} \right| / 2= 0$, i.e.
	$$
\eta = \frac{- 2 d + \sqrt{4 d^2 + 2d | \alpha |^{\frac 1n}  \left| \widetilde c_{j} \right| }}{2 d}.
	$$
 We divide the proof into three cases:
\begin{enumerate}
\item[$i.$] 
If $k' \ge k + \eta k^{2/n}$, we obtain 
	$$
	\left\{
	\begin{array}{l}
\displaystyle \left| \lambda_{j,k}-\lambda_{j',k'} \right| \geq \left| \Im \left( \lambda_{j,k}-\lambda_{j',k'} \right) \right| = |\alpha |^{\frac 1n} \left| \widetilde s_{j} \right| \left[ (k')^{2-{\frac{2}{n}}} - k^{2-{\frac{2}{n}}}\right]\\
	\noalign{\smallskip}
\displaystyle \phantom{\left| \lambda_{j,k}-\lambda_{j',k'} \right| } \geq |\alpha |^{\frac 1n} \left| \widetilde s_{j} \right| \left[  \left( k + \eta k^{2/n} \right)^{2-{\frac{2}{n}}} - k^{2-{\frac{2}{n}}} \right] \geq |\alpha |^{\frac 1n} \left| \widetilde s_{j} \right| \eta \left(2 - \frac 2n \right) k .
	\end{array}
	\right.
	$$
In the previous inequalities we have used~\eqref{xM} with $x= k $ and $M = \eta k^{2/n}$. This proves inequality~\eqref{distint} for $k_1 = 1$.

\item[$ii.$] Let us take $k'\in ( k - \eta k^{2/n} ,k + \eta k^{2/n})$. 
Let us take $k_0$ large enough such that $k - \eta k^{2/n}>0$ for all $k\geq k_0$. If $k \ge k_0$ and using that $n \ge 4$ and
	$$
1 + \frac 2n \le 2 - \frac 2n, \quad \frac 4n \le 2 - \frac 2n ,
	$$
we deduce 
	$$
	\left\{
	\begin{array}{l}
\displaystyle \left| \lambda_{j,k}-\lambda_{j',k'} \right| \ge \left| \Re \left( \lambda_{j,k}-\lambda_{j',k'} \right) \right| \\
	\noalign{\smallskip}
\displaystyle \phantom{ \left| \lambda_{j,k}-\lambda_{j',k'} \right| }\geq  | \alpha |^{\frac 1n}  \left| \widetilde c_{j} \right|  \left[ k^{2-{\frac{2}{n}}}+ (k')^{2-{\frac{2}{n}}} \right] - d \left| k^2-(k')^2  \right| \\
	\noalign{\smallskip}
\displaystyle \phantom{ \left| \lambda_{j,k}-\lambda_{j',k'} \right| } \geq  | \alpha |^{\frac 1n}  \left| \widetilde c_{j} \right| k^{2-{\frac{2}{n}}} -d \max \left\{ k^2 - \left( k - \eta k^{2/n} \right)^2, \left( k + \eta k^{2/n}\right)^2 - k^2 \right\}  \\
	\noalign{\smallskip} 
\displaystyle \phantom{ \left| \lambda_{j,k}-\lambda_{j',k'} \right| } \geq | \alpha |^{\frac 1n}  \left| \widetilde c_{j} \right| k^{2-{\frac{2}{n}}} - 2 d \eta k^{1 + \frac 2n} - d \eta^2 k^{\frac 4n}  \\
	\noalign{\smallskip} 
\displaystyle \phantom{ \left| \lambda_{j,k}-\lambda_{j',k'} \right| }  \ge \left[ | \alpha |^{\frac 1n}  \left| \widetilde c_{j} \right| - 2 d \eta - d \eta^2 \right] k^{2-{\frac{2}{n}}} = \frac 12 | \alpha |^{\frac 1n}  \left| \widetilde c_{j} \right|  k^{2-{\frac{2}{n}}} \,.
	\end{array}
	\right.
	$$
    In the last equality, we have used the expression of $\eta$.
Then, we obtain inequality~\eqref{distint} for $k_1 = k_0 \geq 1$.

\item[$iii.$] Finally, let us consider $k,k'$ such that $1 \le k' \le k - \eta k^{2/n}$. In this case we can repeat the previous arguments. Indeed, if we choose $k \ge k_1 \ge k_0$ large enough, one has:
	$$
	\left\{
	\begin{array}{l}
\displaystyle \left| \lambda_{j,k} - \lambda_{j',k'} \right| \geq \left| \Im \left( \lambda_{j,k} - \lambda_{j',k'} \right) \right| = |\alpha |^{\frac 1n} \left| \widetilde s_{j} \right| \left[ k^{2-{\frac{2}{n}}} - (k')^{2 - {\frac{2}{n}}} \right]\\
	\noalign{\smallskip}
\displaystyle \phantom{\left| \lambda_{j,k}-\lambda_{j',k'} \right| } 
\geq |\alpha |^{\frac 1n} \left| \widetilde s_{j} \right| \left[  k^{2-{\frac{2}{n}}} - \left( k - \eta k^{2/n} \right)^{2-{\frac{2}{n}}}\right] \\
	\noalign{\smallskip}
\displaystyle \phantom{\left| \lambda_{j,k}-\lambda_{j',k'} \right| }
\geq |\alpha |^{\frac 1n} \left| \widetilde s_{j} \right| \eta \left(2 - \frac 2n \right) k^{\frac 2n} \left( k - \eta k^{2/n} \right)^{1 -{\frac{2}{n}}} \ge \widetilde C k,
	\end{array}
	\right.
	$$
where $\widetilde C$ is a new positive constant only depending on $n$, $\alpha$ and $d$. In the previous inequalities we have used that $n \ge 4$ and inequality~\eqref{xM} with $x= k - \eta k^{2/n} $ and $M = \eta k^{2/n}$. This proves inequality~\eqref{distint} for $k_1\geq k_0 \geq 1$.

\end{enumerate}

\smallskip

\noindent\textbf{Case B.4:} $\lambda_{j,k} \in \C \setminus \R$ and $\lambda_{j' ,k' } $ is such that $\widetilde s_j \not= \widetilde s_{j'}$ and $\text{sign} (\widetilde s_j) = \text{sign} (\widetilde s_{j'})$. Observe that in this case, $\widetilde s_j$ $\widetilde s_{j'}$ are non null and we can repeat the same proof as \textbf{Case A.3}. 
This ends the proof. 
\end{proof}

\begin{remark}\label{rmk:1}
	Under the conditions of Proposition~\ref{prop: index 0} we infer that
	$$
\mu(\lambda_{j,k})=1,\footnote{$\mu(\lambda)$ is the {\color{black} geometric} multiplicity of $\lambda\in \sigma(L^\star)$} \quad \forall k \ge k_1, \quad \forall j : 0 \le j \le n-1
	$$
	 and it is easy to see that there exists $\widetilde k_0 \in \N^\star$ such that $\Re \left( \lambda_{j, k} \right) >0$, for any $k \ge \widetilde k_0 $ and any $j \in \{ 0, 1, \dots , n-1 \}$. 
	 Thus, if we define the set 
	\begin{equation}\label{Jset}
\mathcal{J} := \left\{(j,k)\in \{0,\ldots,n-1\}\times \mathbb{N}^\star: \mu(\lambda_{j,k})=1, \quad \Re \left( \lambda_{j, k} \right) >0 \right\},
	\end{equation}
then, under conditions of Proposition~\ref{prop: index 0}, there exists $k_0\in \mathbb{N}^\star$ such that 
	$$ 
\{0,\ldots,n-1\}\times \{k_0, k_0 + 1, \ldots\} \subset \mathcal{J}. 
	$$
On the other hand, the situation is very different when $n$ is even and $\alpha >0$. For instance, when $n = 2$ and $d $ and $\alpha$ does not satisfy~\eqref{f9}, the operator $L^\star$ has an infinite number of eigenvalues with geometric multiplicity equal to $2$ (see~\eqref{f8}).
\end{remark}

\begin{remark}\label{r22}
        
As said in Remark~\ref{r1.4}, when $n = 2p$, $p \in \N^\star$, and $\alpha > 0$, the operator $L^\star :=- D^\star \partial_{xx}f + A^\star $ has, for any $k \ge 1$, exactly two real eigenvalues, $\lambda_{0,k}$ and $\lambda_{p,k}$, and $n-2$ complex eigenvalues (see Proposition~\ref{eigenvectors}). In this case the real eigenvalues are given by
	$$
\lambda_{0,k} = dk^2 + \alpha^{1/n} k^{2 - 1/p} \quad \hbox{and} \quad \lambda_{p ,k} = dk^2 - \alpha^{1/n} k^{2 - 1/p } ,\quad \forall k \ge 1. 
	$$
In this case and in view of the proof of Proposition~\ref{prop: index 0}, we can conclude that inequality~\eqref{cond lambda} holds except for $(j,j')=(0,p)$ or $(j,j')=(p,0)$. 
\end{remark}

As a direct consequence of Proposition~\ref{prop: index 0}, Lemma~\ref{l2} and Remark~\ref{rmk:1}, we have:
\begin{corollary}\label{c1}
Under conditions of Proposition~\ref{prop: index 0}, the sequence $\Lambda = \{\lambda_{j,k}\}_{(j,k)\in \mathcal{J}} $, with $\mathcal{J}$ given by~\eqref{Jset}, satisfies~\eqref{H:lambda} and $c(\Lambda) = 0$.
\end{corollary}

\section{Boundary null controllability}\label{sec:bound}
In this section, we will prove the boundary null controllability results stated in Theorems~\ref{thm:bord} and \ref{tn=2}.

\subsection{Proof of Theorem~\ref{thm:bord}}

The proof of Theorem~\ref{thm:bord} will be developed in three sections. In the first section we will prove the sufficient condition stated in item 1. The second section is devoted to the proof of the necessary condition in item~1. Finally, we will prove item~2 in Section~\ref{Teo1:item2}.

 

\subsubsection{Sufficient condition of the item \texorpdfstring{$1$}{TEXT} of Theorem~\ref{thm:bord}}\label{s4.1}

%
Recall that the matrices $D, A \in \mathcal{L}(\R^n)$ and $ B \in \R^n$ are given by~\eqref{def A B D}, with $d \ge 1$ and $\alpha \in \R^\star$, and the expression of the eigenvalues of $L^\star =- D^\star \partial_{xx} + A^\star $ is given in item a) of Proposition~\ref{eigenvectors}. Remember also that in the case $\alpha >0$ the dimension $n$ of system~\eqref{eq:syst bound} is odd and, therefore,~\eqref{cond lambda} holds.   

Let us first observe that, under assumptions of Theorem~\ref{thm:bord}, one has
	\begin{equation}\label{neq_cond}
\lambda_{j,k}\neq \lambda_{j',k'}, \quad \forall k,k' \in \mathbb{N}^{\star}, \quad \forall j,j'\in\{0,...,n-1\} \hbox{ with } (j, k) \not= (j', k')
	\end{equation}
    and condition~\eqref{f9} holds when $\alpha < 0$ and $ n $ is odd. Then, the goal is to prove that system \eqref{eq:syst bound} is null (resp., approximately) controllable at time $T$.

On the other hand, without loss of generality, we can assume that $\Re \left( \lambda_{j,k} \right) > 0 $ for any $k \ge 1$ and $j : 0 \le j \le n-1$. Indeed, taking into account that $\lim_{k \to \infty} \Re \left( \lambda_{j,k} \right) = \infty$, for any $j : 0 \le j \le n-1$, we can conclude the existence of a positive constant $M > 0$ such that 
	$$
\Re \left( \lambda_{j,k} + M \right) > 0 , \quad \forall k \in \mathbb{N}^{\star}, \quad \forall j \in \{0,...,n-1\}.
	$$
Performing the change $\widetilde z = e^{- M t} z$ in system~\eqref{eq:syst bound}, the controllability properties of this system at time $T$ are equivalent to the corresponding properties of system
\begin{equation*}
	\left\{
	\begin{array}{ll}
		\widetilde z_t  - D\widetilde z_{xx} + (A+MI_n) \widetilde z = 0 & \hbox{in } Q_T,\\
		\noalign{\smallskip}
		\widetilde z(\cdot,0)= B e^{- M t} v, \quad \widetilde z(\cdot,\pi)= 0& \hbox{in } (0,T),\\\noalign{\smallskip}
		\widetilde z(0,\cdot)= z^0						&\hbox{in } (0,\pi),
	\end{array}
	\right.
\end{equation*}
with $z^0\in H^{-1}(0,\pi; \R^n)$ and $v\in L^2(0,T)$. It is clear that 
	$$
\sigma (- D^\star \partial_{xx} + A^\star + M I_n ) = \left\{ \lambda + M : \lambda \in \sigma (- D^\star \partial_{xx} + A^\star ) \right\},
	$$
($I_n \in \mathcal{L} (\R^n)$ is the identity matrix) and then $\Re ( \lambda ) > 0$ for any $\lambda \in \sigma (- D^\star \partial_{xx} + A^\star + M I_n )$.

In order to prove that system~\eqref{eq:syst bound} is null controllable at time $T > 0$, let us first present an equivalent property to the null controllability of the system. Let us introduce the following adjoint system to~\eqref{eq:syst bound}:
	\begin{equation}\label{nondiag_adj}
	\left\{
	\begin{array}{ll}
-\varphi_t  -D^\star\varphi_{xx} + A^\star \varphi=  0 &\hbox{in } Q_T,\\
	\noalign{\smallskip}
\varphi(\cdot,0)= \varphi(\cdot,\pi)=  0 	&\hbox{on } (0,T),\\
	\noalign{\smallskip}
\varphi(T,\cdot)= \varphi^T	&\hbox{in } (0,\pi),
	\end{array}
	\right.
	\end{equation}
where $\varphi^T\in H_0^1(0,\pi; \C^n)$ and matrices $D$, $A$ and $B$ are defined in \eqref{def A B D}. For any $\varphi^T\in H_0^1(0,\pi; \C^n)$ system~\eqref{nondiag_adj} has a unique solution 
	$$
\varphi \in L^2(0,T; H^2(0, \pi; \C^n) \cap H^1_0(0, \pi; \R^n)) \cap \mathcal{C}^0 (0,T; H_0^1(0,\pi; \C^n))
	$$
which depends continuously on the initial data $\varphi^T $. In fact, if $\varphi^T\in H_0^1(0,\pi; \C^n)$, $v\in L^2(0,T; \C)$ and $y^0\in H^{-1}(0,\pi;\C^n)$, the corresponding solutions $z$ and $\varphi$ of systems~\eqref{eq:syst bound} and~\eqref{nondiag_adj} satisfy
	\begin{equation}\label{f10}
\langle y(T,\cdot), \varphi^T \rangle_{H^{-1}, H_0^1} - \langle y^0, \varphi (0,\cdot) \rangle_{H^{-1}, H_0^1} = \int_0^T \left( v (t)  , B^\star D^\star \varphi_x (t,0) \right)_\C \,dt,
	\end{equation}	
where $\langle \cdot , \cdot \rangle_{H^{-1}, H_0^1}$ is the duality pairing between $H^{-1} (0, \pi; \C^n)$ and $H^1_0(0, \pi; \C^n)$.

From identity~\eqref{f10}, the null controllability problem of system~\eqref{eq:syst bound} can be reformulated as a \textit{moment problem}. More precisely, system~\eqref{eq:syst bound} is null controllable at time $T > 0$ if and only if for any initial data $y^0\in H^{-1}(0,\pi;\C^n)$ there exists a control $v \in L^2(0,T ; \C)$ such that 
	\begin{equation*}
- \langle y^0, \varphi (0,\cdot) \rangle_{H^{-1}, H_0^1} = \int_0^T \left( v (t)  , B^\star D^\star \varphi_x (t,0) \right)_\C \,dt, \quad \forall \varphi^T\in H_0^1(0,\pi; \C^n),
	\end{equation*}	
with $\varphi$ the solution of the adjoint problem~\eqref{nondiag_adj} associated to $\varphi^T$. 
Since $\mathcal{B}^\star$ is a Schauder basis of $H^1_0(0,\pi;\mathbb{C}^n)$ (see~\eqref{basis} and Lemma~\ref{riesz_bases}), the null controllability of system~\eqref{eq:syst bound} is equivalent to the following property: 
	
\textbf{Property:} For any initial data $y^0\in H^{-1}(0,\pi;\C^n)$, there exists a control $v\in L^2(0,T ; \C)$ such that 
	\begin{equation}\label{prob_1}
\!\!\!\!\!\!\!-\langle y^0,\varphi_{j,k}(0,\cdot)\rangle_{H^{-1},H^1_0} = \int_0^T \left( v(t) , B^\star D^\star\partial_x\varphi{_{j,k}}(t,0) \right)_\C \,dt,\,~~ \forall k\in\mathbb{N}^\star ,\,~~ \forall j\in \{0,...,n-1\}, 
	\end{equation}	
	where $\varphi_{j,k}$ is the solution of system \eqref{nondiag_adj} associated to the initial data $\varphi^T=\Phi_{j,k} $. 

A simple computation leads to the formula
	\begin{equation*}
\varphi_{j,k}(t,x)=e^{-\lambda_{j,k}(T-t)}\Phi_{j,k} (x),\quad (t,x)\in Q_T,
	\end{equation*}
	whence
\begin{equation*}
	\left\{
	\begin{array}{l}
\displaystyle \varphi_{j,k}(0,x)=e^{-\lambda_{j,k}T}\Phi_{j,k} (x), \quad x\in (0,\pi), \\
	\noalign{\smallskip}
\displaystyle \partial_x\varphi_{j,k}(t,0) = k\sqrt{\frac{2}{\pi}}e^{-\lambda_{j,k}(T-t)} V_{j,k}, \quad \forall t\in (0,T),
	\end{array}
	\right.
	\end{equation*}
where the vector $V_{j,k} \in \C^n$ is given in Proposition~\ref{eigenvectors}. Thus, using these expressions in problem~\eqref{prob_1} we can conclude that system~\eqref{eq:syst bound} is null controllable at time $T > 0$ if and only if 
	\begin{equation}\label{fb}
	\left\{
	\begin{array}{l}
\hbox{there is  }v\in L^2(0,T; \C)  \hbox{ such that for all }k\in\mathbb{N}^\star\hbox{ and } j\in\{0,\ldots, n-1\}, \\
	\noalign{\smallskip}
\displaystyle \left( B^\star D^\star V_{j,k} \right)^\star \int_0^Tv(T-t)e^{-\lambda^\star_{j,k}t}\,dt = - \frac{1}{k}\sqrt{\frac{\pi}{2}} e^{-\lambda^\star_{j,k}T} \langle y^0,\Phi_{j,k} \rangle_{H^{-1},H^1_0} .	
	\end{array}
	\right.
\end{equation}
Let us analyse the expression $B^\star D^\star V_{j,k}$. From~\eqref{def A B D} and Proposition~\ref{eigenvectors}, we deduce that, for each
$k \in \N^\star$ and $j : 0 \le j \le n-1$:
	\begin{equation}\label{eq:dalphaB*}
	\begin{array}{ll}
 \text{$B^\star D^\star V_{j,k} = 0$}
	\end{array}~~ \Longleftrightarrow~~\left\{
	\begin{array}{ll}
\displaystyle 1 + d \alpha^{- \frac 1n} k^{\frac 2n } e^{- \frac{2 \pi j}n i }= 0, & \hbox{if } \alpha >0,  \\
	\noalign{\smallskip}
\displaystyle 1 + d \left| \alpha \right|^{- \frac 1n} k^{\frac 2n } e^{- \frac{(2j +1 ) \pi }n i }= 0, & \hbox{if } \alpha < 0.
	\end{array}
	\right.
	\end{equation}

If $\alpha > 0$, thanks to the assumption $n = 2p + 1$, with $p \in \N^\star$, we can conclude that $B^\star D^\star V_{j,k} \not= 0$ for every $k \in \N^\star$ and $j : 0 \le j \le n-1$. On the other hand, if $\alpha < 0$, $B^\star D^\star V_{j,k} = 0$ if and only if $n  = 2p + 1$, with $p \in \N^\star$, $ j = p$ and $ k =: {\sqrt{\left| \alpha \right|}}/{d^{ n/2}} \in \N^\star$. Thanks to assumption~\eqref{f9} we can also conclude that $B^\star D^\star V_{j,k} \not= 0$ for every $k \in \N^\star$ and $j : 0 \le j \le n-1$. Therefore, problem~\eqref{fb} is equivalent to: 
	\begin{equation}\label{prob_3}
	\left\{
	\begin{array}{l}
\hbox{find }v\in L^2(0,T; \C)  \hbox{ such that for all, }k\in\mathbb{N}^\star\hbox{ and } j\in\{0,\ldots, n-1\}, \\
	\noalign{\smallskip}
\displaystyle\int_0^Tv(T-t)e^{-\lambda^\star_{j,k}t}\,dt = e^{-\lambda^\star_{j,k}T} M_{j,k}(y^0),	
	\end{array}
	\right.
\end{equation}
    where 
$$
    M_{j,k}(y^0):=- \frac{1}{k}\sqrt{\frac{\pi}{2}}\frac{\langle y^0,\Phi_{j,k} \rangle_{H^{-1},H^1_0}}{ V_{j,k}^\star D B } .
$$
	This is the moment problem associated to the boundary null controllability of system~\eqref{eq:syst bound}.
	
	In order to solve the moment problem~\eqref{prob_3}, we will apply Theorem~\ref{lemma: base ortho} to the sequence  $ \Lambda := \left\{ \lambda_{j,k} \right\}_{ k \ge 1, 0 \le j \le n-1 }$ of eigenvalues of the operator $L^\star = - D^\star \partial_{xx}f+ A^\star $. Thanks to condition~\eqref{neq_cond}, the sequence $\Lambda$ satisfies condition~\eqref{H:lambda}. If we use the notation
	\begin{equation}\label{f0'}
p_{j,k}(t):=e^{-\lambda_{j, k}t}, \quad \forall t\in(0,T),
	\end{equation}
we can apply Theorem~\ref{lemma: base ortho} and deduce the existence of a biorthogonal family $\left\{ q_{j,k} \right\}_{ k \ge 1,  0 \le j \le n-1 }$ to $\left\{ p_{j,k } \right\}_{ k \ge 1,  0 \le j \le n-1 }$ in $L^2 (0,T ; \C)$ which satisfies~\eqref{estim q_k}, i.e., a family such that
	\begin{equation*}
\displaystyle \int_0^T p_{j,k} (t) q^\star_{l,m}(t) \, dt = \delta_{km} \delta_{jl},\quad \forall k,m\geq1, \quad \forall j,l :  0 \le j,l \le n-1 .
	\end{equation*}
Moreover, under the assumptions of Theorem~\ref{thm:bord}, we can also apply Proposition~\ref{prop: index 0} and Lem\-ma~\ref{l2} to obtain that $c ( \Lambda ) = 0$. Therefore, the biorthogonal family $\left\{ q_{j,k} \right\}_{ k \ge 1,  0 \le j \le n-1 }$ satisfies the following property: for any $\varepsilon>0$, there exists a constant $C_\varepsilon > 0$ such that
	\begin{equation}\label{estim_q_j,k}
\|q_{j, k}\|_{L^2(0,T; \C)}\leq C_\varepsilon e^{\varepsilon \Re(\lambda_{j, k})}, \quad \forall k \ge 1, \quad \forall j :  0 \le j \le n-1 .
	\end{equation}

We are in conditions to solve problem~\eqref{prob_3}. The function
	$$
v(t)=\sum\limits_{k=1}\limits^{\infty}\sum\limits_{j=0}\limits^{n-1} e^{-\lambda^\star_{j,k}T} M_{j,k}(y^0) q_{j,k}(T-t)
	$$ 
provides a formal solution to this problem. Let us see that, in fact, it is a solution to the moment problem, i.e., let us see that $v \in L^2(0,T; \C)$. 

First, since $y^0 \in H^{-1}(0,\pi;\mathbb{R}^n)$ and taking into account the expression of the vectors $V_{j,k}$ (see Proposition~\ref{eigenvectors}), we infer that, for any $\varepsilon > 0 $, there exists a positive constant $C_\varepsilon$ such that
	$$
\left| M_{j,k}(y^0) \right| \le C_{\varepsilon} e^{\varepsilon \Re(\lambda_{j,k})}\|y^0\|_{H^{-1}(0,\pi;\mathbb{C}^n)}, \quad \forall k \ge 1, \quad \forall j :  0 \le j \le n-1 .
	$$
Let us take $\varepsilon > 0$ (which will be chosen later). Using the previous estimate together with inequality~\eqref{estim_q_j,k}, we get
	\begin{equation*}
	\left\{
	\begin{array}{l}
\displaystyle \|v\|_{L^2(0,T; \C )}\leq \sum_{k=1}^{+\infty}\sum _{j= 0 }^{n- 1 } C_{\varepsilon} e^{-\Re \left( \lambda_{j,k} \right) T}e^{2 \varepsilon \Re(\lambda_{j,k})} \|y^0\|_{H^{-1}(0,\pi;\mathbb{C}^n)} \\
	\noalign{\smallskip}
\displaystyle \phantom{\|v\|_{L^2(0,T; \C )}}
= \sum_{k=1}^{+\infty}\sum _{j = 0}^{n-1} C_{\varepsilon} e^{-\Re \left( \lambda_{j,k} \right) \left( T - 2 \varepsilon \right)} \|y^0\|_{H^{-1}(0,\pi;\mathbb{C}^n)} .
	\end{array}
	\right.
	\end{equation*}
where $C_\varepsilon$ is a positive constant. Taking, for example, $\varepsilon = T/4$, we obtain that the series in the definition of $v$ converges absolutely in $L^2 (0, T ; \C)$. Thus, the previous control $v$ solves the moment problem~\eqref{prob_3}. This proves the null controllability result at time $T>0$ of system~\eqref{eq:syst bound}.


\subsubsection{Necessary condition of Item 1 of Theorem ~\ref{thm:bord}}\label{s4.2}


	
 First of all, notice that the approximate controllability at time $T$ for system \eqref{eq:syst bound} is equivalent to a  Fattorini-Hautus test. More precisely,
\begin{theorem}\label{theo:fattorini}
System~\eqref{eq:syst bound} is approximately controllable at time $T$ if and only if, for every $\lambda\in\mathbb{C}$ and $\Phi\in D( L^\star)$, we have the following property
	\begin{equation*}
	\left. 
	\begin{array}{ll}
L^\star\Phi  = \lambda\Phi &\hbox{in } (0,\pi)\\
	\noalign{\smallskip}
B^\star D^\star \partial_x \Phi(0)= 0 &
	\end{array}\right\}
\Longrightarrow \Phi = 0 \hbox{ in } (0,\pi).
	\end{equation*}
	\end{theorem}
\noindent For the proof, one just have to apply \cite[Theorem $1.1$]{O14}.

	Let us use Theorem~\ref{theo:fattorini} applied to the operator $L = - D \partial_{xx} + A$.
	To do this, by contradiction, assume first that $n = 2 p +1$, with $p \in \N^\star$, $\alpha < 0$ and ${\sqrt{\left| \alpha \right|}}/{d^{ n/2}} := \mathcal{K}  \in \N^\star$. In this case, $ \Phi = \Phi_{p,\mathcal{K}} $ (see item \textit{a)} of Proposition~\ref{eigenvectors}) satisfies $\Phi \not \equiv 0$, $L^\star\Phi  = \lambda\Phi $, with $\lambda = \lambda_{p,\mathcal{K}} $, and 
	$$
B^\star D^\star \partial_x \Phi(0) =  \mathcal{K} \sqrt{\frac{2}{\pi}} B^\star D^\star V_{p , \mathcal{K}} = - \mathcal{K} \sqrt{\frac{2}{\pi}} | \alpha |^{- \frac{n-2}n}k^{2 \frac{n-2}n }\left( 1 - d \left| \alpha \right|^{- \frac 1n} \mathcal{K}^{\frac 2n } \right)   = 0 .
	$$
From Theorem~\ref{theo:fattorini}, we deduce that system~\eqref{eq:syst bound} is not approximately controllable at time $T$.

	On the other hand, assume now that there exists $k,k'\in\mathbb{N}^{\star}$ and $j,j'\in\{0,...,n-1\}$ such that $(j,k) \not= (j', k')$ and 
	$$
\lambda:= \lambda_{j,k} = \lambda_{j', k'}.
	$$
	The eigenvalue $\lambda$ is, at least, double (geometric multiplicity) and some associated eigenfunctions are 
	$$
\widetilde \Phi_{j,k}=\widetilde V_{j,k}w_{k}\quad\hbox{and}\quad\widetilde \Phi_{j',k'}=\widetilde V_{j',k'}w_{k'},
	$$
where $\widetilde{V}_{j,k}:=\left( c_{j,k}^l \right)_{1 \le l \le n} \in \C^n$ and the coefficients $c_{j,k}^l$, $1 \le l \le n$, are given in~\eqref{expr phi_k} (see~Proposition~\ref{eigenvectors}).

	Let us point out that  $B^\star D^\star=(0, \cdots, 0,1, d)$. Thus, we introduce 
	$$
\Phi:=k' \left[c_{j',k'}^{n-1} + dc_{j',k'}^n \right] \widetilde \Phi_{j,k} - k \left[c_{j,k}^{n-1} + d c_{j,k}^n \right]\widetilde \Phi_{j',k'}.
	$$
	It is not difficult to see that $\Phi$ is not identically zero and satisfies
	\begin{equation*}
	\left\{ 
	\begin{array}{ll}
-D^\star\partial_{xx}  \Phi+A^\star\Phi  = \lambda\Phi &\hbox{in } (0,\pi)\\
	\noalign{\smallskip}
B^\star D^\star\partial_x\Phi(0)= 0.
	\end{array}
	\right.
	\end{equation*}
Therefore, Theorem~\ref{theo:fattorini} leads to the non-approximate controllability of system~\eqref{eq:syst bound}. 
\hfill $\Box$


\subsubsection{Proof of the item $2$ of Theorem ~\ref{thm:bord}}\label{Teo1:item2}


%
As in Section \ref{s4.1}, we will apply the moment method to prove the null controllability at time $T>0$ of system~\eqref{eq:syst bound} when $\alpha = 0$. {In particular, $L^\star = - D^\star \partial_{xx}$ and 
from Proposition ~\ref{eigenvectors}, we have}
	$$
\sigma ( L^\star ) = \{ \lambda_k := dk^2 : k \in \N^\star \}.
	$$

As in the previous case, system~\eqref{eq:syst bound} is null controllable at time $T > 0$ if and only if there exists a control $v \in L^2(0,T)$ such that~\eqref{prob_1} holds, where $\Phi_{j,k} = e_{n - j} w_k$ (see item b) of Proposition~\ref{eigenvectors}).

For the initial data $\varphi^T:=\Phi_{j,k}  = e_{n - j} w_k$, $k \ge 1$, $j: 0 \le j \le n-1$, the solution to the adjoint problem~\eqref{nondiag_adj} is given by:
    	\begin{equation}\label{explicity_adj_sol jordan}
\varphi_{j,k}(t,x) = e^{-dk^2 (T-t)} \sum_{l = 0}^{j} (-1)^{j - l} \frac{k^{2(j - l)}}{(j - l)!} (T-t)^{j - l} e_{n - l} w_k (x), \quad (t,x) \in Q_T.
	\end{equation}
From this identity, we infer
	\begin{equation*}
	\left\{
	\begin{array}{l}
\displaystyle \varphi_{j,k}(0,x) = e^{-d k^2 T}\sum_{l = 0}^{j} (-1)^{j - l} \frac{\left( k^2T \right)^{ j - l}}{(j - l)!} e_{n - l} w_k (x) ,\\
	\noalign{\smallskip}
\displaystyle \partial_x \varphi_{j,k}(t,0) = k\sqrt{\dfrac{2}{\pi}} e^{-dk^2 (T-t)} \sum_{l = 0}^{j} (-1)^{j - l} \frac{k^{2(j - l)}}{(j - l)!} (T-t)^{j - l} e_{n - l} ,
	\end{array}
	\right.
	\end{equation*}
for any $k \ge 1$ and $j: 0 \le j \le n-1$.

Thus, if we introduce the notation
	\begin{equation}\label{f1}
\widetilde M_{j,k}(y^0) = - \frac{1}{ k}\sqrt{\frac{\pi}{2}} \sum_{l = 0}^{j} (-1)^{ l} \frac{\left( k^2T \right)^{ j - l}}{(j - l)!} \langle y^0 , e_{n - l} w_k \rangle_{H^{-1},H^1_0}, \quad \forall k \in \N^\star, \ 0 \le j \le n -1
	\end{equation}
and observing that
	$$
B^\star D^\star = \left( 0, 0, \dots, 1, d \right),
	$$
the moment problem~\eqref{prob_1} becomes: find $v\in L^2(0,T)$ such that for all $k \in \N^\star$ one has
	$$
	\left\{
	\begin{array}{l}
\displaystyle d \int_0^T v(T-t) e^{-d k^2 t} \,dt = e^{-d k^2 T} \widetilde M_{0 ,k}(y^0) , \\
	\noalign{\smallskip}
\displaystyle - \frac{k^{2(j - 1)}}{(j - 1) !} \int_0^T v(T-t) t^{j - 1} e^{-d k^2 t} \,dt + d \frac{k^{2j}}{j !} \int_0^T v(T-t) t^j e^{-d k^2 t} \,dt = e^{-d k^2 T} \widetilde M_{j,k}(y^0), \\
	\noalign{\smallskip}
\displaystyle \forall j :  1 \le j \le n - 1.
	\end{array}
	\right.
	$$
This linear system can be written in \MD{a} vectorial form as
	\begin{equation}\label{f2}
\mathcal{A}_k X_k =  e^{-d k^2 T}  \widetilde M_{k}(y^0),
	\end{equation}

where
	\begin{equation}\label{f3}
	\begin{array}{c}
\mathcal{A}_k = \left(
	\begin{array}{cccccc}
d & 0 & 0 & 0 & \cdots & 0  \\
-1 & dk^2 & 0 & 0 & \cdots & 0 \\
0 & - k^2 & \frac d2 k^4 & 0 & \cdots & 0 \\
0 & 0 & - \frac 12 k^4 & \frac d6 k^6& \cdots  & 0 \\
\vdots & \vdots & \vdots & \ddots & \ddots & \vdots \\
0 & 0 & 0 & \cdots & - \frac{k^{2(n - 2)}}{(n - 2) !} & d \frac{k^{2(n - 1)}}{(n - 1) !}
	\end{array}
\right) \in \mathcal{L} (\R^n), \\
	\noalign{\smallskip} 
\displaystyle X_k = \left( \int_0^T v(T-t) t^j e^{-d k^2 t} \,dt \right)_{0 \le j \le n-1} \in \R^n, \quad  \widetilde M_{k}(y^0)= \left( \widetilde M_{j,k}(y^0) \right)_{0 \le j \le n-1} \in \R^n.
	\end{array}
	\end{equation}

System~\eqref{f2} is triangular, then it is equivalent to 
	\begin{equation}\label{prob_3 jordan}
	\left\{
	\begin{array}{l}
\hbox{Find }v\in L^2(0,T)  \hbox{ such that }\\
	\noalign{\smallskip}
\displaystyle \int_0^T v(T-t) t^{j}e^{-\lambda_{k}t} \,dt =  e^{-d k^2 T} M_{j,k}(y^0) , \quad \forall k \in \N^\star, \ 0 \le j \le n -1,
	\end{array}
	\right.
\end{equation}
where the coefficients $M_{j,k}(y^0)$ are given by 
	\begin{equation}\label{f4}
\left( M_{j,k}(y^0) \right)_{0 \le j \le n-1} := \mathcal{A}_k^{-1} \widetilde M_{k}(y^0) \in \R^n .
	\end{equation}
\MD{In} summarizing, the null controllability property of system~\eqref{eq:syst bound} at time $T$ is equivalent to the moment problem~\eqref{prob_3 jordan}. 

Our next step will be to prove that the moment problem~\eqref{prob_3 jordan} admits a solution $v \in L^2 (0, T)$. Firstly, we can apply Theorem~\ref{t3}, with $\eta = n$, to the real sequence $\Lambda := \{ d k^2 \}_{k \ge 1}$ and deduce the existence of a family $\left\{ q_{j,k} \right\}_{k\geq 1, 0 \leq j \leq n - 1}\subset L^2 (0,T)$ biorthogonal to $\left\{t^{j}e^{-\lambda_{k}t} \right\}_{k\geq 1, 0 \leq j \leq n - 1}$ which satisfies~\eqref{biacot}. As in the previous case, this fact provides a formal solution to the moment problem~\eqref{prob_3 jordan}:
	$$
v(t) := \sum_{k=1}^{\infty} \sum_{j=0}^{n-1} e^{-d k^2 T} M_{j,k}(y^0) q_{j,k}(T-t).
	$$ 

On the other hand, the previous series converges absolutely in $L^2 (0,T)$. Indeed, taking into account the expression of the coefficients $M_{j,k}(y^0)$ (see~\eqref{f4},~\eqref{f1} and~\eqref{f2}), the estimate~\eqref{biacot} and $y^0 \in H^{-1}(0, \pi ; \R^n )$, it is not difficult to prove the following property: for any $\varepsilon > 0$, there exists a positive constant $C_\varepsilon$ such that 
	$$
\left| M_{j,k}(y^0) \right| \| q_{j,k} \|_{L^2(0, T )} \le C_\varepsilon e^{\varepsilon d k^2 } \|y^0\|_{H^{-1}(0,\pi;\R^n)}, \quad \forall (j, k) : k \ge 1, \ 0 \le j \le n - 1.
	$$
With this inequality, we can reason as in Section \ref{s4.1} and prove that $v$ is a solution to the moment problem~\eqref{prob_3 jordan}. 
\hfill $\Box$


\subsection{Proof of Theorem~\ref{tn=2}}\label{s4}
We will devote this section to prove Theorem~\ref{tn=2}. To this end, let us consider system~\eqref{eq:syst bound} in the case $n = 2$, with matrices $D, A \in \mathcal{L}(\R^2)$ and $ B \in \R^2$ given by~\eqref{def A B D}, with $d \ge 1$ and $\alpha > 0 $. In this case, recall Proposition~\ref{eigenvectors}, the eigenvalues of the operator $L^\star = - D^\star \partial_{xx} + A^\star $ are given by
	\begin{equation}\label{f7}
\lambda_{0,k}:=d k^2 + \sqrt \alpha k , \quad \lambda_{1,k}:=d k^2 - \sqrt \alpha k, \quad k \ge 1 
	\end{equation}
{and $\Phi_{0,k} = V_{0,k} w_k $ and $\Phi_{1,k} = V_{1,k} w_k $, with $ V_{j,k} = \widetilde V_{j,k}/ |\widetilde V_{j,k} |$ ($j =0,1$) and
	\begin{equation}\label{f11}
\widetilde V_{0,k} = \left( 
	\begin{array}{c}
1 \\ \alpha^{-1/2} k
	\end{array}
	\right), \quad 
\widetilde V_{1,k} = \left( 
	\begin{array}{c}
1 \\ - \alpha^{-1/2} k
	\end{array}
	\right),
	\end{equation}
are eigenvectors associated to $\lambda_{0,k}$ and $\lambda_{1,k}$.}

As in Section~\ref{s4.1} and without loss of generality, we are going to assume that $\sigma (L^\star) \subset (0, \infty)$, i.e., $\lambda_{1,k} > 0$ for all $k \in \N^\star$.

Observe that in the case $n = 2$, inequalities~\eqref{cond lambda0} and~\eqref{H:gap} are, in general, not valid when one takes as sequence $\Lambda$ the real sequence $\Lambda:= \left\{ \lambda_{0,k}, \lambda_{1,k} \right\}_{k \ge 1}$. Indeed, from~\eqref{f7}, we deduce the equalities
	\begin{equation}\label{f8}
	\left\{
	\begin{array}{l}
\displaystyle \lambda_{1, k + m} - \lambda_{0,k} = \left( 2 k + m\right) \left( dm - \sqrt \alpha \right) , \quad \forall k, m \in \N^\star , \\
	\noalign{\smallskip}
\displaystyle \lambda_{1, k} - \lambda_{1 ,\ell } = \left(  k - \ell \right) \left[ d(k + \ell)  - \sqrt \alpha \right] , \quad \forall k, \ell \in \N^\star.
	\end{array}
	\right.
	\end{equation}
Thus, if ${\sqrt \alpha}/d = m \in \N^\star$, we have $\lambda_{0,k} = \lambda_{1, k +m} $, for any $k \ge 1$, and $\lambda_{1, k} = \lambda_{1,  m - k}$, for any $k : 1 \le k \le m-1$.

As a consequence, we deduce that the operator $L^\star$ has an infinite number of eigenvalues with geometric multiplicity equal to $2$. Therefore, we can follow the arguments of Section~\ref{s4.2} and conclude that system~\eqref{eq:syst bound} is neither approximately nor null controllable in $H^{-1}(0, \pi ; \R^2)$ at any time $T>0$. This proves the necessary part of Theorem~\ref{tn=2}.

\smallskip 

Let us now assume that ${\sqrt \alpha}/d \not\in \N^\star$ and prove that system~\eqref{eq:syst bound} is null controllable in the space $H^{-1} (0, \pi; \R^2)$ at time $T>0$. As in Section~\ref{s4.1}, this controllability result is equivalent to the moment problem:
	\begin{equation}\label{f12}
	\left\{
	\begin{array}{l}
\hbox{Given $y^0 \in H^{-1} (0, \pi; \R^2)$, find }v\in L^2(0,T)  \hbox{ such that for all }k\in\mathbb{N}^\star\hbox{ and } j = 0,1, \\
	\noalign{\smallskip}
\displaystyle B^\star D^\star V_{j,k} \int_0^Tv(T-t)e^{-\lambda_{j,k}t}\,dt = - \frac{1}{k}\sqrt{\frac{\pi}{2}} e^{-\lambda_{j,k}T} \langle y^0,\Phi_{j,k} \rangle_{H^{-1},H^1_0} .	
	\end{array}
	\right.
\end{equation}
We follow the arguments of Section~\ref{s4.1} in order to solve the previous moment problem.

In this case, we can write
	$$
m-1 < \frac{\sqrt \alpha}d < m,
	$$
where $m \ge 1$ is a positive integer. From~\eqref{f8}, it is not difficult to prove that the eigenvalues of $L^\star$ are simple and one has
	$$
\lambda_{1, k+m-1} < \lambda_{0, k} < \lambda_{1, k + m}, \quad \forall k \ge 1,
	$$
and
	$$
	\left\{
	\begin{array}{l}
\lambda_{1, k + m} - \lambda_{0, k} = \left( 2 k + m\right) \left( dm - \sqrt \alpha \right) \ge 3 \left( dm - \sqrt \alpha \right) := \mathcal C_0 > 0, \\
	\noalign{\smallskip}
\lambda_{0, k} - \lambda_{1, k + m -1}  = \left( 2 k + m - 1 \right) \left[ \sqrt \alpha - d (m-1)\right] \ge 2 \left[ \sqrt \alpha - d (m-1)\right] := \mathcal C_1 > 0,
	\end{array}
	\right.
	$$
for any positive integer $ k \in \N^\star$. Thus, the sequence $\Lambda = \left\{ \lambda_{0,k}, \lambda_{1,k} \right\}_{k \ge 1}$ can be rearranged as an increasing sequence of positive real numbers $\Lambda = \left\{ \Lambda_{k} \right\}_{k \ge 1}$ as follows:
	$$
	\left\{
	\begin{array}{ll}
\left\{ \Lambda_{k} \right\}_{1 \le k \le m} := \left\{ \lambda_{1, k} \right\}_{1 \le k \le m} , \\
	\noalign{\smallskip}
\Lambda_{m + 2\ell -1 } := \lambda_{0, \ell }, & \forall \ell \ge 1, \\
	\noalign{\smallskip}
\Lambda_{m + 2\ell } := \lambda_{1, m + \ell }, & \forall \ell \ge 1 . 
	\end{array}
	\right.
	$$

Therefore, the increasing sequence $\Lambda$ satisfies properties~\eqref{H:lambda} and~\eqref{H:gap} for 
	$$
\rho := \min \{ \mathcal C_0, \mathcal C_1,  \mathcal C_2\} > 0, \quad  \mathcal C_2 := \frac 1m \min_{1 \le k, \ell \le m } \left| \lambda_{1,k} - \lambda_{1,\ell} \right|.
	$$ 
On the other hand, as in Section~\ref{s4.1} (see\eqref{eq:dalphaB*}), we can see that the assumption ${\sqrt \alpha}/d \not\in \N^\star$ implies that $B^\star D^\star V_{j,k} \not= 0$ for every $k \in \N^\star$ and $j = 0, 1$ (see~\eqref{f11}). Therefore, we can apply the arguments of Section~\ref{s4.1} to solve the moment method~\eqref{f12} and obtain the null controllability property of this system in $H^{-1}( 0 , \pi ; \R^2)$ at time $T > 0$. This proves the sufficient part of Theorem~\ref{tn=2}. 

\smallskip

Let us now prove the last part of Theorem~\ref{tn=2}. So, assume that 
	\begin{equation}\label{f13}
\frac{\sqrt \alpha}d = m \in \N^\star .
	\end{equation}
As said before, the operator $L^\star$ has an infinite number of eigenvalues with geometric multiplicity equal to $2$ (see~\eqref{f8}) and, in general, the moment problem~\eqref{f12} cannot be solved. Nevertheless, this problem has a solution $v \in L^2(0,T)$ for some initial data $y^0 \in H^{-1}(0, \pi; \R^2)$. Let us see this point.

Firstly,  we have $\lambda_{0,k} = \lambda_{1, k +m} $, for any $k \ge 1$, and $\lambda_{1, k} = \lambda_{1,  m - k}$, for any $k : 1 \le k \le m-1$. Thus, the real sequence $\Lambda:= \left\{ \lambda_{0,k}, \lambda_{1,k} \right\}_{k \ge 1}$ is, in fact, $\Lambda = \left\{ \lambda_{1,k} \right\}_{k \ge m_0} $, where $m_0 = 1 + \left\lfloor \frac{m-1}2 \right\rfloor$ ($\left\lfloor \cdot \right\rfloor$ is the floor function). From the expression of $m$ (see~\eqref{f13}) and $m_0$, we can prove that the sequence $\Lambda $ is increasing and satisfies properties~\eqref{H:lambda} and~\eqref{H:gap}.

Secondly, let us remember that, given $y^0 \in H^{-1}(0, \pi; \R^2)$, there exists $v \in L^2(0,T)$ such that the solution $y$ of~\eqref{eq:syst bound} satisfies $y (T, \cdot) = 0$ in $(0, \pi)$ if and only if the control $v$ solves the moment problem~\eqref{f12}. Using the expressions~\eqref{f11} and~\eqref{f13}, this moment problem can be rewritten as 
	\begin{equation}\label{f14}
	\left\{
	\begin{array}{l}
\hbox{Given $y^0 \in H^{-1} (0, \pi; \R^2)$, find }v\in L^2(0,T)  \hbox{ such that for all }k\in\mathbb{N}^\star , \\
	\noalign{\smallskip}
\displaystyle \left( 1 + \frac km\right) \int_0^Tv(T-t)e^{-\lambda_{0,k}t}\,dt = - \frac{1}{k}\sqrt{\frac{\pi}{2}} e^{-\lambda_{0,k}T} \left( y^0_k, \widetilde V_{0,k} \right )_{\R^2} ,	\\
	\noalign{\smallskip}
\displaystyle \left( 1 - \frac km\right) \int_0^Tv(T-t) e^{-\lambda_{1,k} t} \, dt = - \frac{1}{k}\sqrt{\frac{\pi}{2}} e^{-\lambda_{1 ,k}T} \left( y^0_k, \widetilde V_{1 ,k} \right )_{\R^2} ,
	\end{array}
	\right.
	\end{equation}
where $y_k^0 = \langle y^0, w_k \rangle_{H^{-1},H^1_0} \in \R^2$, $k \ge 1$, are the Fourier coefficients of $y^0$.

For $y^0 \in H^{-1} (0, \pi ; \R^2)$, let us consider the conditions
	\begin{equation}\label{f15}
	\left\{
	\begin{array}{ll}
\displaystyle \left( y_k^0, \widetilde V_{1,k} \right)_{\R^2} = \left( y_{m-k}^0, \widetilde V_{1,m-k} \right)_{\R^2} & \forall k : 1 \le k \le m_0 -1,\\
	\noalign{\smallskip}
\displaystyle \left( y_k^0, \widetilde V_{0,k} \right)_{\R^2} = - \left( y_{k + m}^0, \widetilde V_{1, k + m} \right)_{\R^2} & \forall k \ge 1,\\
	\noalign{\smallskip}
\displaystyle \left( y_m^0, \widetilde V_{1,m} \right)_{\R^2} = 0, &
	\end{array}
	\right.
	\end{equation}
and introduce the closed subspace of $H^{-1}(0, \pi; \R^2)$ given by:
	$$
\mathcal{X} := \left\{ y^0 \in H^{-1}(0, \pi; \R^2) : y^0 \hbox{ satisfies conditions }\eqref{f15} \right\} .
	$$

The set $\mathcal{X}$ is a closed subspace of $H^{-1}(0, \pi; \R^2)$ which has infinite codimension. Indeed, if we consider the closed subspace
	$$
\mathcal{Y} := \left\{ e_1 f : f_k = 0,~ k \in \{1, \dots, m \}, \ \left\{f_k \right\}_{k \ge 1} \hbox{ are the Fourier coefficients of } f \in H^{-1}(0, \pi)\right\}
	$$
then, thanks to conditions~\eqref{f15}, one has
	$$
\mathcal{X} \cap \mathcal{Y} = \{ 0 \}.
	$$
On the other hand, taking into account that $\dim \mathcal{Y} = \infty$, we deduce that $\mathcal{X}$ has infinite codimension\footnote{$H^{-1}(0, \pi; \R^2)$ is a Hilbert space. Then, consider $P : H^{-1}(0, \pi; \R^2) \to \mathcal{X}$, the orthogonal projection of $H^{-1}(0, \pi; \R^2) $ onto $\mathcal{X}$. Thus, $\mathcal A : y \in \mathcal{Y} \mapsto \mathcal A ( y ) = y - Py \in \mathcal{X}^\bot$ is injective.}.
%
%

In order to finish the proof of Theorem~\ref{tn=2}, let us see the property: \textit{``given $y^0 \in H^{-1}(0, \pi; \R^2)$, there exists a control $ v \in L^2 (0,T)$ such that the solution $y $ of system~\eqref{eq:syst bound} satisfies $y ( T, \cdot) = 0$ in $(0, \pi)$ if and only if $y^0 \in \mathcal X$''}, i.e., let us prove that the moment problem~\eqref{f14} has solution if and only if $y^0 \in \mathcal X$.

Let us start assuming that $y^0 \in \mathcal{X}$. In this case, from~\eqref{f15} and the expressions of $m$ and the functions $\Phi_{j,k}$ (see~\eqref{f13} and~\eqref{f11}), the moment problem~\eqref{f14} is equivalent to 
	\begin{equation*}
	\left\{
	\begin{array}{l}
\hbox{Find }v\in L^2(0,T)  \hbox{ such that for all }k \ge m_0 , \ k \not= m , \\
	\noalign{\smallskip}
\displaystyle \left( 1 - \frac km\right) \int_0^Tv(T-t) e^{-\lambda_{1,k} t} \, dt = - \frac{1}{k}\sqrt{\frac{\pi}{2}} e^{-\lambda_{1 ,k}T} \left( y^0_k, \widetilde V_{1 ,k} \right )_{\R^2} .
	\end{array}
	\right.
	\end{equation*}

As proved before, the  sequence $\Lambda = \left\{ \lambda_{1,k} \right\}_{k \ge m_0} $ is increasing and satisfies properties~\eqref{H:lambda} and~\eqref{H:gap}. Following the arguments of Section~\ref{s4.1}, we deduce that the previous moment problem admits a solution $v \in L^2(0,T)$.

Let us now suppose that $y^0 \not\in \mathcal{X}$. In this case, the moment problem is incompatible and does not admit any solution:

\begin{enumerate}
\item If $\left( y_m^0, \widetilde V_{1,m} \right)_{\R^2} \not= 0$, it is clear that the second equation of~\eqref{f14} has no solution when $k = m$.
\item If for some $k_0 : 1 \le k_0 \le m_0 - 1$ (resp., $k_0 \ge 1 $) one has $\left( y_{k_0}^0, \widetilde V_{1,k_0} \right)_{\R^2} \not= \left( y_{m-k_0}^0, \widetilde V_{1,m-k_0} \right)_{\R^2}$ (resp., $\left( y_{k_0}^0, \widetilde V_{0,k_0} \right)_{\R^2} \not= - \left( y_{k_0 + m}^0, \widetilde V_{1, k_0 + m} \right)_{\R^2}$), from the equality $\lambda_{1, k_0} = \lambda_{1,  m - k_0}$ (resp., $\lambda_{0,k_0} = \lambda_{1, k_0 +m}$), it is not difficult to show that the problem~\eqref{f14} is incompatible.  
\end{enumerate}
This proves the previous equivalence and ends the proof of Theorem~\ref{tn=2}. \hfill $\Box$


\subsection{Pointwise controllability}\label{s4.3}
In this section we will prove the null controllability at time $T>0$ of system~\eqref{f16} in $L^2(0, \pi; \R^n)$ when $x_0 \in (0, \pi)$ satisfies appropriate properties. To this end, we will follow the same ideas of the proof of Theorem~\ref{thm:bord}. 
%

\begin{remark}
Taking into account that $\delta_{x_0} \in H^{-1}(0, \pi)$, we deduce that, for any $y^0 \in L^2 (0, \pi; \R^n)$ and $u \in L^2(0,T)$, system~\eqref{f16} admits a unique solution $y$ with regularity
	$$
y\in L^2(0,T;H_0^1(0,\pi;\R^n))\cap  \mathcal{C}^0([0,T]; L^2(0,\pi;\R^n)),
	$$
and which depends continuously on $y^0$ and $u$. 
\end{remark}

Let us first describe the approximate controllability result for system~\eqref{f16}. One has:

%
%

\begin{theorem}\label{pointwise-app}
Let us consider the matrices $D, A \in \mathcal{L}(\R^n)$ and $ B \in \R^n$ given by~\eqref{def A B D}, with $d \ge 1$ and $\alpha \in \R$. In addition, assume that $n = 2p+1$ ($p\in\mathbb{N}^{\star}$) when $\alpha >0$. Let us also fix $x_0 \in (0, \pi)$. Then, 

\begin{enumerate}
\item If $\alpha \not= 0$, system~\eqref{f16} is approximately controllable at time $T$ if and only if  one has:
	\begin{equation}\label{f17}
	\left\{
	\begin{array}{l}
\hbox{the eigenvalues of } L^\star=- D^\star \partial_{xx} + A^\star \hbox{ has geometric multiplicity equal to }1;  \\
	\noalign{\smallskip}
x_0 \not= r \pi, \hbox{ with }r \in \Q \cap (0,1).
	\end{array}
	\right.
	\end{equation} 
\item If $\alpha =0$, system~\eqref{f16} is approximately controllable at time $T > 0$ if and only if $x_0 \not= r \pi$, with $r \in \Q \cap (0,1) $.
\end{enumerate}
\end{theorem}

\begin{proof}
As saw in Section~\ref{s4.2}, it is not to difficult to show that system~\eqref{f16} is approximately controllable in $L^2(0, \pi; \R^2)$ at time $T > 0$ if and only if, for every $\lambda\in\C$ and $\Phi\in D( L^\star)$, we have the following property
	\begin{equation}\label{f18}
	\left. 
	\begin{array}{ll}
L^\star\Phi  = \lambda\Phi &\hbox{in } (0,\pi)\\
	\noalign{\smallskip}
B^\star \Phi(x_0) = 0 &
	\end{array}\right\}
\Longrightarrow \Phi = 0 \hbox{ in } (0,\pi).
	\end{equation}
We will do the proof for $\alpha \not= 0 $. The case $\alpha = 0 $ can be obtained following the same argument.

Let us first see that the conditions in~\eqref{f17} are necessary for the approximate controllability of system~\eqref{f16}. We can argue as in Section~\ref{s4.2} and prove that the first condition in~\eqref{f17} is necessary. On the other hand, if $x_0 = r \pi$, with $ r = m / \ell \in \Q \cap (0,1)$ and $m, \ell \in \N^\star$, we deduce that $\lambda = \lambda_{0, \ell} \in \R$ and $\Phi := \Phi_{0, \ell} = V_{0, \ell} w_\ell $ (see Proposition~\ref{eigenvectors}) satisfy $L^\star \Phi = \lambda \Phi$, 
	$$
B^\star \Phi (x_0) = \sqrt{\frac{2}{\pi}} B^\star V_{0, \ell} \sin ( \ell x_0  ) = \sqrt{\frac{2}{\pi}} B^\star D^\star V_{0, \ell} \sin (m \pi ) = 0,
	$$
and $\Phi \not\equiv 0$ in $(0, \pi)$. Therefore, from~\eqref{f18} we deduce that system~\eqref{f16} cannot be approximately controllable at time $T > 0$.

\smallskip

Let us see that conditions in~\eqref{f17} imply the approximate controllability of system~\eqref{f16} at time $T>0$ i.e., the Fattorini-Hautus test~\eqref{f18}. Indeed, first, if $\Phi\in D( L^\star)$ satisfies $L^\star\Phi = \lambda \Phi$ in $(0, \pi)$, for $\lambda \in \C$, we deduce $\lambda = \lambda_{j,k}$ and $\Phi = C \Phi_{j, k}$, with $k \ge 1$, $0 \le j \le n-1$ and $C \in \C$. Secondly, taking into account the expression of $\Phi_{j, k}$ (see Proposition~\ref{eigenvectors}),  condition $B^\star \Phi(x_0) = 0$ can be also written  as
    $$
0 = B^\star \Phi (x_0) =  \left\{
	\begin{array}{ll}
\displaystyle \frac{C}{\left| \left( c_{j,k}^l \right)_{1 \le l \le n} \right|} \sqrt{ \frac{2}{\pi}} \left( \alpha^{-\frac{n -1}n} k^{2{\frac{(l-1)}{ n}}}e^{- \frac{2\pi j}n (n -1)i} \right) \sin ( k x_0  ) , & \hbox{if } \alpha >0,  \\
	\noalign{\smallskip}
\displaystyle  \frac{C}{\left| \left( c_{j,k}^l \right)_{1 \le l \le n} \right|} \sqrt{ \frac{2}{\pi}} \left( | \alpha |^{-\frac{n -1}n} k^{2{\frac{( n -1)}{ n}}}e^{- \frac{(2 j + 1) \pi}n (n -1)i}  \right) \sin ( k x_0  ) , & \hbox{if } \alpha < 0.
	\end{array}
	\right.
	$$
In both cases, assumptions~\eqref{f17} imply $C = 0$ and, therefore, $\Phi \equiv 0$. Then, one has~\eqref{f18} and the approximate controllability of system~\eqref{f16} at time $T>0$. This ends the proof.
\end{proof}

%
%

%
%
\begin{remark}
Taking into account Proposition~\ref{prop: index 0} (see~Remark~\ref{rmk:1}), we deduce that, when $\alpha \not= 0$ is under the assumptions of Theorem~\ref{pointwise-app}, the eigenvalues of $L^\star$ satisfy $\mu(\lambda_{j,k})=1$ for any $k \ge k_1$ and any $j : 0 \le j \le n - 1$. Thus, if $x_0 \in (0, \pi)$ is such that $x_0/ \pi \not \in \Q $, system~\eqref{f16} is approximately controllable in $L^2(0, \pi; \R^2)$ at time $T > 0$ apart from a finite dimensional space of $L^2(0, \pi; \R^2)$.

On the other hand, the proof of Theorem~\ref{pointwise-app} is still valid when $\alpha >0$ and $n = 2p$, with $p\in\mathbb{N}^{\star}$, that is to say, system~\eqref{f16} is approximately controllable at time $T$ if and only if~\eqref{f17} holds. Nevertheless, in this case, we can have an infinite number of eigenvalues $\lambda$ of $L^\star $ such that $\mu (\lambda ) \ge 2$ (this is the case when $n=2$ and $p$ and $\alpha$ satisfy~\eqref{f13}, see Section~\ref{s4}). So, we cannot conclude that, if $x_0 \in (0, \pi)$ is such that $x_0/ \pi \not \in \Q $, system~\eqref{f16} is approximately controllable in $L^2(0, \pi; \R^2)$ at time $T > 0$ apart from a finite dimensional space of $L^2(0, \pi; \R^2)$.
\end{remark}

Let us now study the null controllability of system~\eqref{f16}. The result reads as follows:

%
%

%
\begin{theorem}\label{pointwise-null}
Under conditions of Theorem~\ref{pointwise-app}, let us assume that~\eqref{f17} holds when $\alpha \not= 0$. In addition, assume that $x_0 = \vartheta \pi$, with $\vartheta \in (0, 1 )$ an irrational number, and consider 
	$$
T_\vartheta = \limsup_{k\to\infty} \frac{ - \log \left\vert \sin ( k \vartheta \pi ) \right\vert }{ d k^2 } \in [0, \infty].
	$$
Then:
\begin{enumerate}
\item System~\eqref{f16} is null controllable in $ L^{2} ( 0,\pi ; \C^n ) $ at any time $T > T_\vartheta $.
\item System~\eqref{f16} is not null controllable in $ L^{2} ( 0,\pi ; \C^n ) $ for $T < T_\vartheta $.
\end{enumerate}
\end{theorem}

\begin{proof} 
Again, we will do the proof of the result when $\alpha \not=0 $. The case $\alpha = 0$ can be obtained from a similar argument (see Section~\ref{Teo1:item2}). So, assume that~\eqref{f17} holds and take $x_0 = \vartheta \pi$ where $\vartheta \in (0, 1 )$ is an irrational number. 

Let us first prove the first item in Theorem~\ref{pointwise-null}. To this end, assume $T_\vartheta \in [0, \infty)$ and take $T > T_\vartheta$. As in Section~\ref{s4.1}, the null controllability at time $T > 0$ of system~\eqref{f16} is equivalent to the following moment problem: given $y^0 \in L^2(0, \pi ;\C^n)$, find a control $v \in L^2(0, \pi; \C )$ such that
%
	\begin{equation*}
\displaystyle \sqrt{\frac{2}{\pi}} \sin ( k \vartheta \pi ) \left( B^\star  V_{j,k} \right)^\star \int_0^Tv(T-t)e^{-\lambda^\star_{j,k}t}\,dt = - e^{-\lambda^\star_{j,k}T} ( y^0,\Phi_{j,k} )_{L^2} , 
	\end{equation*}
for any $k \in \N^\star$ and any $j: 0 \le j \le n-1$. From~ \eqref{f17} we deduce that $\mu (\lambda_{j,k}) = 1$ and $\sin ( k \vartheta \pi )  \not= 0$. Also, it is easy to check that $ B^\star V_{j,k} \not= 0$. Thus, the previous moment problem is equivalent to
	\begin{equation}\label{f19}
	\left\{
	\begin{array}{l}
\hbox{Given }y^0 \in L^2(0, \pi ;\C^n) , \hbox{ find } v \in L^2(0, \pi; \C )  \hbox{ s.t.~for all }k \ge 1 \hbox{ and } j: 0 \le j \le  n-1 , \\
	\noalign{\smallskip}
\displaystyle\int_0^Tv(T-t)e^{-\lambda^\star_{j,k}t}\,dt = e^{-\lambda^\star_{j,k}T} \widetilde M_{j,k}(y^0),	
	\end{array}
	\right.
	\end{equation}
where
$$
 \widetilde M_{j,k}(y^0)  :=- \frac{1}{\sin ( k \vartheta \pi )}\sqrt{\frac{\pi}{2}}\frac{( y^0,\Phi_{j,k} )_{L^2} }{ V_{j,k}^\star B }. 
$$

Again, a formal solution \MD{to} the previous moment problem is 
	$$
v(t) = \sum\limits_{k=1}\limits^{\infty}\sum\limits_{j=0}\limits^{n-1} e^{-\lambda^\star_{j,k}T} \widetilde M_{j,k}(y^0) q_{j,k}(T-t), 
	$$ 
where $\left\{ q_{j,k} \right\}_{ k \ge 1,  0 \le j \le n-1 }$ is a biorthogonal family  to $\left\{ p_{j,k} \right\}_{ k \ge 1,  0 \le j \le n-1 }$ (see~\eqref{f0'}) in $L^2 (0,T ; \C)$ satisfying~\eqref{estim_q_j,k} (see Theorem~\ref{lemma: base ortho}). 

Let us check that the previous series is absolutely convergent in $L^2 (0, T; \C)$. Indeed, from Proposition~\ref{eigenvectors}, we get
	$$
d k^2 - | \alpha |^{1/n} k^{2 - \frac 2n} \le \Re (\lambda^\star_{j,k} ) \le d k^2 + | \alpha |^{1/n} k^{2 - \frac 2n} , \quad \forall k \ge 1, \ j: 0 \le j \le  n-1,
	$$
On the other hand, using the expression of $T_\vartheta$, we deduce that for any $\varepsilon > 0$ there exists a positive constant $C_\varepsilon$ such that
	$$
\dfrac{1}{ | \sin ( k \vartheta \pi ) |} \le C_\varepsilon e^{( T_\vartheta + \varepsilon)dk^2  }, \quad \forall k \ge 1. 
	$$
Finally, repeating the arguments in Section~\ref{s4.1}, we also have that, for any $\varepsilon > 0$ there exists a positive constant $C_\varepsilon$ such that
	$$
\left|\frac{( y^0,\Phi_{j,k} )_{L^2}}{  V_{j,k}^\star B }\right|  \le C_\varepsilon e^{ \varepsilon dk^2  } \|y^0\|_{L^{2}(0,\pi;\mathbb{C}^n)} , \quad \forall k \ge 1 , \ j: 0 \le j \le  n-1.
	$$

Therefore, since $\left\{ q_{j,k} \right\}_{ k \ge 1,  0 \le j \le n-1 }$ satisfies ~\eqref{estim_q_j,k}, we have
	\begin{equation*}
	\left\{
	\begin{array}{l}
\displaystyle \|v\|_{L^2(0,T; \C )}\leq \sum_{k=1}^{\infty}\sum _{j=0}^{n-1} C_{\varepsilon} e^{-\Re \left( \lambda_{j,k} \right) T} e^{( T_\vartheta + \varepsilon)dk^2  } e^{ \varepsilon dk^2  } \|y^0\|_{L^{2}(0,\pi;\mathbb{C}^n)}e^{\varepsilon\Re \left( \lambda_{j,k} \right) } \\
	\noalign{\smallskip}
\displaystyle \phantom{\|v\|_{L^2(0,T; \C )}}
\leq \sum_{k=1}^{\infty}\sum _{j=0}^{n-1} C_{\varepsilon} e^{- d k^2 T} e^{(T+\varepsilon) | \alpha |^{1/n} k^{2 - \frac 2n} }e^{( T_\vartheta + \varepsilon)dk^2  } e^{ 2\varepsilon dk^2  } \|y^0\|_{L^{2}(0,\pi;\mathbb{C}^n)} \\
	\noalign{\smallskip}
\displaystyle \phantom{\|v\|_{L^2(0,T; \C )}}
= \sum_{k=1}^{\infty}\sum _{j=1}^{n} C_{\varepsilon} e^{-dk^2\left( T - T_\vartheta - 3 \varepsilon \right)} e^{(T+\varepsilon) | \alpha |^{1/n} k^{2 - \frac 2n} } \|y^0\|_{L^{2}(0,\pi;\mathbb{C}^n)} .
	\end{array}
	\right.
	\end{equation*}
where $C_\varepsilon$ is a new positive constant. It is clear that, taking $\varepsilon = \left( T - T_\vartheta  \right)/6$, the previous series converges absolutely. This proves that $v \in L^2 (0, T; \C)$ and we have constructed a solution of the moment problem~\eqref{f19}. This shows the first item in Theorem~\ref{pointwise-null}. 

\smallskip

Let us now assume that $ T_\vartheta \in (0, \infty]$ and consider $0 < T < T_\vartheta$. The objective is to prove that system~\eqref{f16} is not null controllable in $ L^{2} ( 0,\pi ; \C^n ) $ at time  $T $. To this end, we will use the following result:
\begin{theorem}\label{t-observab}
Let us consider the matrices $D, A \in \mathcal{L}(\R^n)$ and $ B \in \R^n$ given by~\eqref{def A B D}, with $d \ge 1$ and $\alpha \in \R$. Then, system~\eqref{f16} is not null controllable in $ L^{2} ( 0,\pi ; \C^n ) $ at time  $T > 0$ if and only if there exists a constant $C_T > 0$ such that
	\begin{equation}\label{observab}
\| \varphi(0, \cdot) \|^2_{L^{2} ( 0,\pi ; \C^n ) } \le C_T \int_0^T \left| B^\star \varphi(t, x_0) \right|^2 \, dt , \quad \forall \varphi^T \in L^{2} ( 0,\pi ; \C^n ),
	\end{equation}
where $\varphi \in L^2(0,T ; H_0^1(0, \pi; \C^n)) \cap \mathcal{C}^0([0, T]; L^{2} ( 0,\pi ; \C^n ) )$ is the solution of the adjoint problem~\eqref{nondiag_adj} associated to $\varphi^T \in L^{2} ( 0,\pi ; \C^n )$. 
\end{theorem}

For a proof of this result, see~\cite{T-W} and~\cite{Z}.

Let us see that the observability inequality~\eqref{observab} fails when $T < T_\vartheta$ and, therefore, system~\eqref{f16} is not null controllable in $ L^{2} ( 0,\pi ; \C^n ) $ at time  $T $. By contradiction, assume that, for a positive constant $C_T$, inequality~\eqref{observab} holds. In particular, if we take 
	$$
\varphi^T (x) = V_{0,k} w_k(x) = \sqrt{\frac{2}{\pi}} V_{0,k} \sin(kx), \quad x \in (0, \pi),
	$$
(see Proposition~\ref{eigenvectors}), the corresponding solution of system~\eqref{f16} is given by
	$$
\varphi (t, x) = \sqrt{\frac{2}{\pi}} e^{-\lambda_{0,k}(T-t)} V_{0,k} \sin(kx), \quad (t,x) \in Q_T,
	$$
and inequality~\eqref{observab} becomes
	\begin{equation}\label{observab'}
\left| V_{0,k}  \right|^2 e^{- 2 \Re(\lambda_{0,k}) T } \le {\frac{2}{\pi}}  C_T \left| B^\star V_{0,k}  \right|^2 \frac{1}{2 \Re(\lambda_{0,k})} \left( 1 - e^{- 2 \Re(\lambda_{0,k}) T }\right) \left| \sin(k \vartheta \pi) \right|^2, \quad \forall k \ge 1. 
	\end{equation}

On the other hand, let us take $\varepsilon > 0$ such that $T + \varepsilon < T_\vartheta$. From the definition of $T_\vartheta$ we deduce the existence of a subsequence $\left\{ k_n \right\}_{n \ge 1}$ such that
	$$
T + \varepsilon < \frac{ - \log \left\vert \sin ( k_n \vartheta \pi ) \right\vert }{ d k_n^2 }, \quad \forall n \ge 1,
	$$
or, equivalently,
	$$
\left| \sin(k_n \vartheta \pi) \right|^2 < e^{ - 2 d k_n^2 (T + \varepsilon )}  , \quad \forall n \ge 1.
	$$
Combining the previous inequality and~\eqref{observab'} written for the subsequence $\left\{ k_n \right\}_{n \ge 1}$, we get
	$$
\left| V_{0,k_n}  \right|^2 e^{- 2 \Re(\lambda_{0,k_n }) T } \le {\frac{2}{\pi}}  C_T \left| B^\star V_{0,k_n }  \right|^2 \frac{1}{2 \Re(\lambda_{0,k_n })} \left( 1 - e^{- 2 \Re(\lambda_{0,k_n }) T }\right) e^{ - 2 d k_n^2 (T + \varepsilon )}, \quad \forall n \ge 1,
	$$
that is to say,
	$$
\displaystyle 0 < \frac \pi 2 C_T^{-1 } \le  \mathcal{A}_n , \quad \forall n \ge 1, 
	$$
	where
$$
    \mathcal{A}_n:=\frac{\left| B^\star V_{0,k_n }  \right|^2}{ \left| V_{0,k_n}  \right|^2 }  \frac{1}{2 \Re(\lambda_{0,k_n })} \left( 1 - e^{- 2 \Re(\lambda_{0,k_n }) T }\right) e^{ - 2 d k_n^2 \varepsilon  + 2 [\Re(\lambda_{0,k_n })-d k_n^2] T }   \quad \forall n \ge 1.
$$
Finally, taking into account the expressions of $V_{0,k }$ and $\lambda_{0,k }$ (see Proposition~\ref{eigenvectors}) and the inequality
	$$
\Re (\lambda_{0,k } ) \le dk^2 + |\alpha|^{1/n} k^{2-{\frac{2}{n}}}, \quad \forall k \ge 1,
	$$
we obtain that $\displaystyle\lim_{n \to +\infty}  \mathcal{A}_n = 0$. 

This provides a contradiction and the proof of item 2 of Theorem~\ref{pointwise-null}. 
\end{proof}
\begin{remark}
To the authors' knowledge, the first pointwise null controllability result for a parabolic equation was proved in ~\cite{dolecki} for the one dimensional heat equation. Similar results were obtained for coupled parabolic systems in~\cite{AKBGBdT14}.
\end{remark}


\section{Parabolic systems with non-diagonalizable diffusion matrix and non-constant coefficients}\label{sec_negative_}

This section is devoted to prove Theorem \ref{theo:negatif}. The negative part relies on the Fattorini-Hautus test applied to the operator $L_0 := -D \partial_{xx} + q A_0 $. On the other hand, the positive part relies on the \textit{algebraic resolvability} (see~\cite{gromov}).

\begin{proof}[Proof of Theorem~\ref{theo:negatif}, item $a)$]
    First of all, notice that, from Theorem~\ref{theo:fattorini}, sys\-tem \eqref{syst simpl} is approximately controllable at time $T$ if and only if the following property for the adjoint operator $L_0^\star$ holds:
	\begin{equation}\label{eq:non_approx}
	\begin{array}{l}
	\text{
 For every $\lambda \in \mathbb{C}$ and $(\psi, \varphi)\in H^2(0,\pi;\mathbb{R}^2)\cap  H^1_0(0,\pi;\mathbb{R}^2)$, it holds}\\
 	\noalign{\smallskip}
	\left. 
	\begin{array}{ll}
L_0^\star \left( \begin{array}{l}
     \psi  \\
     \varphi 
\end{array} \right) = \lambda \left( \begin{array}{l}
     \psi  \\
     \varphi 
\end{array} \right) & \quad\hbox{in}\quad (0,\pi)\\
	\noalign{\smallskip}
B^\star \left( \begin{array}{l}
     \psi  \\
     \varphi 
\end{array} \right) = 0 &\quad\hbox{in}\quad\omega
	\end{array}
	\right\}
\Longrightarrow \left( \begin{array}{l}
     \psi  \\
     \varphi 
\end{array} \right) = \left( \begin{array}{l}
     0  \\
     0 
\end{array} \right) \quad\text{in}\quad (0,\pi).
	\end{array}
	\end{equation}
    The idea is to construct a potential $q \in \mathcal{C}^{\infty}([0,\pi]) $ such that property~\eqref{eq:non_approx} does not holds. 
    To this end, consider $\omega:=(5\pi/12,7\pi/12)$ and let us construct three functions $\varphi,~\psi\in H^2(0,\pi)\cap H_0^1(0,\pi)$ and $q\in\mathcal{C}^{\infty}( [0, \pi])$ satisfying
	\begin{equation}\label{cond:constr simpl}
	\left\{ 
	\begin{array}{ll}
-\partial_{xx}\psi + q(\cdot) \psi= 36\psi&\hbox{in } (0,\pi),\\
	\noalign{\smallskip}
-\partial_{xx}  \varphi-\partial_{xx} \psi  = 36\varphi &\hbox{in } (0,\pi),\\
	\noalign{\smallskip}
\varphi= 0& \hbox{in }\omega,\\
	\noalign{\smallskip}
\psi \not \equiv 0, \ \varphi \not \equiv 0 &\hbox{in }(0,\pi).
	\end{array}
	\right.
	\end{equation}
%
The strategy will be to construct a suitable function $\psi$ as a perturbation of $\sin(2x)$. With such a function, we define $q$ and $\varphi$ as \eqref{cond:constr simpl}$_1$ and~\eqref{cond:constr simpl}$_2$, respectively, and we check that $(\varphi,\psi)$ satisfies~\eqref{cond:constr simpl}. 
    
    Consider $\psi$ a function of $\mathcal{C}^{\infty}([0,\pi])\cap H_0^1(0,\pi)$ satisfying
	\begin{equation}\label{cond:constr simpl2}
	\left\{
	\begin{array}{ll}
\displaystyle \psi(x) = \sin(2x) + C_1 \theta_1 (x) - C_2 \theta_2 (x) & \displaystyle \forall x \in \left[ 0, \frac{\pi}{3} \right] \cup \left[ \frac{2\pi}{3}, \pi \right], \\
	\noalign{\smallskip}
\displaystyle \psi(x) = - \frac{6}{\pi} x+3 & \displaystyle \forall x \in \overline{\omega} =\left[\frac{5 \pi}{12}, \frac{7\pi}{12}\right],\\
	\noalign{\smallskip}
\displaystyle | \psi(x) - \sin(2x)| < \varepsilon & \displaystyle \forall x \in \left[ \frac{\pi}{3} , \frac{5 \pi}{12} \right]\cup \left[ \frac{7 \pi}{12} , \frac{2\pi}{3} \right],
	\end{array}
	\right.
	\end{equation}
where  $\theta_1$ and $\theta_2$ are two nontrivial nonnegative functions of $\mathcal{C}^{\infty}([0,\pi])$ satisfying
	\begin{equation}\label{fsupp}
	\left\{ 
	\begin{array}{l}
\displaystyle \text{supp}(\theta_1)\subset \left(\frac{\pi}{24}, \frac{3\pi}{24} \right),\\
	\noalign{\smallskip}
\displaystyle \text{supp}(\theta_2)\subset \left( \frac{21\pi}{24} , \frac{23\pi}{24} \right),
	\end{array}
	\right.
	\end{equation}
$\varepsilon>0$ is small enough and $C_1$ and $C_2$ are two positive constants to be determined. The graph of $\psi$ is given in Figure~\ref{fig:contre exemple}.
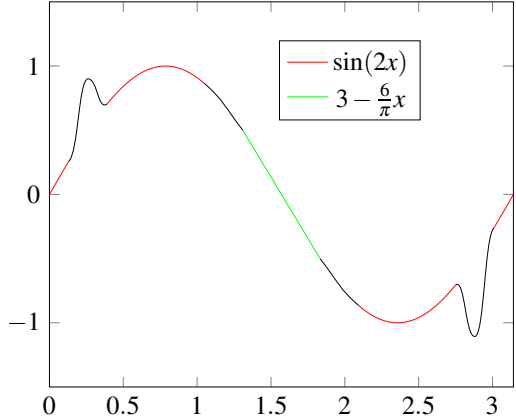
\begin{figure}[H]\begin{center}
    \begin{tikzpicture}[scale=0.9]
\begin{axis}[xmin=0,xmax=pi,ymin=-1.5,ymax=1.5,legend style={at={(0.8,0.9)}}]
\addplot[samples=500,domain=0:pi/24,red]{sin(deg(2*x))};
\addplot[samples=500,domain=5*pi/12:7*pi/12,green]{3-6*x/pi};
\addplot[samples=500,domain=8*pi/12:21*pi/24,red]{sin(deg(2*x))};
\addplot[samples=500,domain=23*pi/24:pi,red]{sin(deg(2*x))};
\addplot[samples=500,domain=3*pi/24:4*pi/12,red]{sin(deg(2*x))};
\draw (axis cs: pi/24,{sin(2*pi/24 r)}) .. controls (axis cs: 3*pi/48,{sin(2*pi/24 r)+pi/48}) and (axis cs: 3*pi/48,{sin(2*2*pi/24 r)+0.4}) .. (axis cs: 2*pi/24,{sin(2*2*pi/24 r)+0.4});
\draw  (axis cs: 2*pi/24,{sin(2*2*pi/24 r)+0.4}) .. controls (axis cs: 5*pi/48,{sin(2*2*pi/24 r)+0.4}) and (axis cs: 5*pi/48,{sin(2*3*pi/24 r)-pi/48}) .. (axis cs: 3*pi/24,{sin(2*3*pi/24 r)});
\draw  (axis cs: 4*pi/12,{sin(2*4*pi/24 r)}) .. controls (axis cs: 9*pi/24,{sin(2*4*pi/24 r)-pi/24}) and (axis cs: 9*pi/24,0.5+1.5*pi/24) .. (axis cs: 5*pi/12,0.5);
\draw  (axis cs: 7*pi/12,-0.5) .. controls (axis cs: 15*pi/24,-0.5-1.5*pi/24) and (axis cs: 15*pi/24,-0.7351257) .. (axis cs: 8*pi/12,{sin(2*8*pi/12 r)});
\draw  (axis cs: 21*pi/24,{sin(2*21*pi/24 r)}) .. controls (axis cs: 43*pi/48,{sin(2*21*pi/24 r)+pi/48}) and (axis cs: 43*pi/48,{sin(2*21*pi/24 r)-0.4}) .. (axis cs: 22*pi/24,{sin(2*21*pi/24 r)-0.4});
\draw  (axis cs: 22*pi/24,{sin(2*21*pi/24 r)-0.4}) .. controls (axis cs: 45*pi/48,{sin(2*21*pi/24 r)-0.4}) and (axis cs: 45*pi/48,{sin(2*23*pi/24 r)-pi/48}) .. (axis cs: 23*pi/24,{sin(2*23*pi/24 r)});
 \legend{$\sin(2x)$,$3- \frac{6}{\pi} x$}
\end{axis}
\end{tikzpicture}
\caption{Example of a function $\psi$ in $[0,\pi]$}\label{fig:contre exemple}
\end{center}
\end{figure}

Notice that the function $\varphi\in\mathcal{C}^{\infty}([0,\pi])$, defined by
	\begin{equation*}
\varphi(x) := \left[\tau  - \frac16 \int_0^x \cos(6y) \psi_{yy}(y)\, dy \right] \sin (6x) + \left[ \frac16 \int_0^x \sin(6y)\psi_{yy}(y) \, dy \right] \cos(6x),  \quad \forall x\in [0,\pi],
	\end{equation*}
satisfies the second equation of \eqref{cond:constr simpl}, where $\tau\in\mathbb{R}$ is a constant to be fixed later.

Let us now verify the boundary conditions and \eqref{cond:constr simpl}$_3$ for $\varphi$.
	Let us first prove that $C_1$, in~\eqref{cond:constr simpl2}, and $\tau$ can be chosen such that $\varphi\equiv 0$ in $\omega$.
	As $\psi(x)=- \frac6\pi x+3$ for $x\in\omega$ and $\psi$ coincides with $\sin(2x)$ in a neighborhood of $0$, we have 
	\begin{equation}\label{f6}
	\left\{
	\begin{array}{l}
\displaystyle - \frac16 \int_0^{\frac{5\pi}{12}} \!\!\! \cos(6y) \psi_{yy}(y)\, dy = 6 \int_0^{\frac{5\pi}{12}} \!\!\! \cos(6y)\psi(y)\, dy - \frac 16, \\ 
	\noalign{\smallskip}
\displaystyle \frac16 \int_0^{\frac{5\pi}{12}} \!\!\! \sin(6y)\psi_{yy}(y) \, dy = -6 \int_0^{\frac{5\pi}{12}}\!\!\!\sin(6y)\psi(y)dy - \frac 1\pi,
	\end{array}
	\right.
	\end{equation}
and for all $x\in\omega$:
	\begin{equation*}
\varphi(x) = \left[\tau- \frac16 + 6 \int_0^{\frac{5\pi}{12}} \!\!\! \cos(6y)\psi(y)\, dy \right] \sin(6x) - \left[ \frac 1\pi + 6 \int_0^{\frac{5\pi}{12}} \!\!\! \sin(6y)\psi(y) \, dy\right] \cos(6x).
	\end{equation*}

	Since
	$$
\displaystyle \frac{1}{\pi}+6 \int_0^{\frac{5\pi}{12}} \!\!\! \sin(6y)\sin(2y) \, dy = \frac{1}{\pi}- \frac{3\sqrt{3}}{16} <0,
	$$
thanks to \eqref{cond:constr simpl2}$_3$ (for $\varepsilon$ small enough), one can choose $C_1$ (recalling that $\sin(6x)>0$ in the interval $(\pi/24,3\pi/24)$) such that
	\begin{equation}\label{coef cos nul}
\frac{1}{ \pi}+ 6\int_0^{\frac{5\pi}{12}} \!\!\! \sin(6y)\psi(y) \, dy=0.
	\end{equation}
In this way, for 
	$$
\displaystyle \tau: = \frac16 - 6 \int_0^{\frac{5\pi}{12}} \!\!\! \cos(6y)\psi(y)\, dy ,
	$$
we obtain $\varphi=0$ in $\omega$.

Let us now verify the boundary conditions for $\varphi$. Notice that $\varphi(0)=0$, by definition.  A suitable choice of $C_2$ will give us $\varphi(\pi)=0$. Indeed, using that $\psi$ is an affine function in $\omega$ and from equalities~\eqref{f6} and~\eqref{coef cos nul}, we have
	$$
	\begin{array}{l}
\displaystyle \varphi ( \pi ) = \frac 16 \int_0^\pi \sin(6y)\psi_{yy}(y) \, dy = \frac16 \int_0^{\frac{5\pi}{12}} \!\!\! \sin(6y)\psi_{yy}(y) \, dy + { \frac1 6}\int_{\frac{7\pi}{12}}^\pi \sin(6y) \psi_{yy}(y) \, dy  \\
	\noalign{\smallskip}
\displaystyle \phantom{\varphi(\pi)} = { \frac1 6}\int_{\frac{7\pi}{12}}^\pi \sin(6y) \psi_{yy}(y) \, dy=-{\frac 1 \pi }-6 \int_{\frac{7\pi}{12}}^\pi\sin(6y)\psi(y) \, dy.
	\end{array}
	$$

Since
	$$
-{\frac 1 \pi} - 6\displaystyle \int_{\frac{7\pi}{12}}^\pi \sin(6y) \sin(2y)\, dy= - \frac 1 \pi + \frac{3\sqrt{3}}{16}>0
	$$
again thanks to \eqref{cond:constr simpl2}$_3$ (for $\varepsilon$ small enough) one can choose $C_2$ (recalling that $\sin(6x)<0$ in $(21\pi/24,23\pi/24)$) such that
	$$
- \frac 1 \pi -6 \displaystyle \int_{\frac{7\pi}{12}}^\pi \sin(6y)\psi(y) \, dy=0
	$$
and then $\varphi(\pi)=0$.

Finally, to verify the first equality in \eqref{cond:constr simpl}, we define $q\in\mathcal{C}^{\infty}([0, \pi])$ by
	\begin{equation}\label{fq}
q := \dfrac{\partial_{xx} \psi + 36\psi}{\psi},
	\end{equation}
with $\psi$ given in~\eqref{cond:constr simpl2}. Taking into account that $\psi$ is null only at points $0$, $\pi/2$ and $\pi$, we infer the existence of neighborhoods of $0$ and $\pi$ in which $\psi$ is equal to $\sin(2x)$ and a neighborhood of $\pi/2$ in which $\psi$ is equal to $- { \frac 6\pi}x+3$. Therefore, we have that $q$ is equal to $32$ in the neighborhoods of $0$ and $\pi$ and equal to $36$ in the neighborhood of $\pi/2$. 
Therefore, the function $q$ is bounded and item $a)$ in Theorem~\ref{theo:negatif} is proved. 
\end{proof}


\begin{remark}
The previous proof provides an argument to construct functions $\varphi,~\psi\in H^2(0,\pi)\cap H_0^1(0,\pi)$ and $q\in\mathcal{C}^{\infty}( [0, \pi])$ satisfying~\eqref{cond:constr simpl}. In fact, this construction is valid for any functions $\theta_1, \theta_2 \in \mathcal{C}^{\infty}([0,\pi])$ satisfying~\eqref{fsupp}. We will use this construction in the proof of  Theorem~\ref{theo:negatif}, item $b)$. 
\end{remark}


Let us now prove item b) of Theorem \ref{theo:negatif}. To do that, we will use the following result, whose proof is given below, after the proof of this theorem.

\begin{theorem}\label{th:fic}
Let $T>0$, $\omega\subset (0,\pi)$, a nonempty open set, and $q\in\mathcal{C}^{\infty}([0,\pi])$ such that $q$ is not constant on an open subset $\omega_1 \subset \omega$. Then, system~\eqref{syst simpl} is null controllable (then, approximately controllable)  at time $T$.
\end{theorem}

\begin{proof}[Proof of Theorem~\ref{theo:negatif}, item $b)$]
Let us take $\omega = (\pi/24,3\pi/24)$ and $\theta_1, \theta_2 \in \mathcal{C}^{\infty}([0,\pi])$ satisfying~\eqref{fsupp} and $\theta_1 (x) =e^x$ for any $x \in \omega_1$, with $\omega_1$ an open interval such that $\omega_1 \subset \! \subset(\pi/24,3\pi/24) $. With the previous choice, let us consider the function $q\in\mathcal{C}^{\infty}([0, \pi])$ given in~\eqref{fq}, with $\psi$ given in~\eqref{cond:constr simpl2}. It is clear that the functions $\varphi,~\psi\in H^2(0,\pi)\cap H_0^1(0,\pi)$ and $q\in\mathcal{C}^{\infty}( [0, \pi])$ satisfy~\eqref{cond:constr simpl} and  
	$$
q(x) =37 - 5\dfrac{\sin(2x)}{\sin(2x)+C_1e^x} , \quad \forall x\in\omega_1.
	$$
Since $q$ is not constant on $\omega_1$, we can apply Theorem \ref{th:fic} and guarantee that system~\eqref{syst simpl} is null controllable at time $T>0$. This ends the proof of item $b)$.
\end{proof}

Let us now return to Theorem~\ref{th:fic} and establish its proof.

\begin{proof}[Proof of Theorem \ref{th:fic}]
Since $q$ is not constant on the open set $\omega_1 \subset\!\subset \omega$, there exist a constant $C >0$ and a new open subset $\widehat{\omega} \subset\omega_1$ such that 
	\begin{equation}\label{eq:partialxq}
|\partial_xq| > C > 0 \mbox{ and } |q| > C > 0 \mbox{ on }\widehat{\omega}.
	\end{equation}
Now, let us reduce the proof of the null controllability of system~\eqref{syst simpl} is null controllable at time $T$ to the resolution of two problems:

\begin{itemize}
\item \textit{Analytic problem:} find $(\widehat{y},\widehat{u})$, with $\widehat{y} \in L^2(0,T;H_0^1(0,\pi;\R^2)) \cap \mathcal{C}^0([0,T]; L^2(0,\pi;\mathbb{R}^2))) $ and $ \widehat{u} \in \mathcal{C}^k(Q_T;\mathbb{R}^2)$ ($k$ is a positive integer that will be fixed later) such that: 
	\begin{equation}\label{eq:ana prob}
	\left\{
	\begin{array}{ll}
\partial_t\widehat{y}_1 = \partial_{xx}\widehat{y}_1 + \partial_{xx}\widehat{y}_2 - q(\cdot) \widehat{y}_1 + \widehat{u}_1 & \hbox{in } Q_T,\\
	\noalign{\smallskip}
\partial_t\widehat{y}_2=\partial_{xx}\widehat{y}_2 + \widehat{u}_2 &\hbox{in } Q_T,\\
	\noalign{\smallskip}
\widehat{y}=0&\hbox{on } (0,T)\times\{0,\pi\},\\
	\noalign{\smallskip}
\widehat{y}(0,\cdot)=y^0, \ y(T,\cdot)=0&\hbox{in } (0,\pi),\\
	\noalign{\smallskip}
\text{supp}(\widehat{u})\subset (0,T)\times\widehat\omega,\\
	\noalign{\smallskip}
\widehat{\omega}\subset\!\subset\omega.
	\end{array}
	\right.
	\end{equation}
\item \textit{Algebraic problem:} find $(z,v)\!\in L^2(0,T;H_0^1(0,\pi;\mathbb{R}^2))\cap  \mathcal{C}^0([0,T]; L^2(0,\pi;\R^2 ) ) ) \times L^2(Q_T)$ such that:
	\begin{equation}\label{eq:alg prob}
	\left\{
	\begin{array}{ll}
\partial_t z_ 1= \partial_{xx} z_1 + \partial_{xx} z_2 - q(\cdot) z_1 + \widehat{u}_1 &\hbox{in } (0,T)\times\omega,\\
	\noalign{\smallskip}
\partial_t z_2 =\partial_{xx} z_2 + v + \widehat{u}_2 &\hbox{in } (0,T)\times\omega,\\
	\noalign{\smallskip}
\text{supp}(z,v)\subset \! \subset (0,T)\times\omega.
	\end{array}
	\right.
	\end{equation}
\end{itemize}

If we are able to solve the analytic and algebraic problems,  then $(\widehat{y}-z,-v)$ is a solution to the null controllability problem for system~\eqref{syst simpl}.

The next task will be to solve the analytic and algebraic problems. 

The resolution of the analytic problem~\eqref{eq:ana prob} is standard and can be established taking into account that $q \in \mathcal{C}^{\infty}([0,\pi])$ and thanks to the local regularity of parabolic equations (see~\cite{bodart} and~\cite{GB-PG} where the local regularity is used to construct regular controls).

\medskip

Now, let us present a resolution of the algebraic problem. System \eqref{eq:alg prob} can be rewritten as follow
	$$
\mathcal{L}(z,v)=\widehat{u},
	$$
where
	$$
\mathcal{L}(z,v) = \left(
	\begin{array}{c}
\partial_t z_1 - \partial_{xx}z_1 - \partial_{xx}z_2 + q(\cdot) z_1 \\
	\noalign{\smallskip}
\partial_t z_2 - \partial_{xx} z_2 - v
	\end{array}
	\right).
	$$
Let us search a differential operator $\mathcal{M}$ with $\mathcal{C}^{\infty}$ coefficients such that
	\begin{equation}\label{eq:LM = I}
\mathcal{L}\circ\mathcal{M} = Id,
	\end{equation}
thus $(z,v):=\mathcal{M}(\widehat{u})$ will be a solution to \eqref{eq:alg prob}.
In the analytic problem \eqref{eq:ana prob}, we search for $\widehat{u}$ regular enough ($k$ large enough)
in order to apply the differential operator $\mathcal{M}$.
The formal adjoint to~\eqref{eq:LM = I} is given by
	\begin{equation}\label{eq:L*M*}
\left( \mathcal{M}^{\star}\circ\mathcal{L}^\star \right) \psi = \psi,
	\end{equation}
where 
	$$
\mathcal{L}^\star(\psi) = \left( 
	\begin{array}{c}
\mathcal{L}_1^\star\psi\\
	\noalign{\smallskip}
\mathcal{L}_2^\star\psi\\
	\noalign{\smallskip}
\mathcal{L}_3^\star\psi
	\end{array}
	\right)
= \left(
	\begin{array}{c}
-\partial_t \psi_1 - \partial_{xx}\psi_1 + q (\cdot) \psi_1 \\
	\noalign{\smallskip}
-\partial_t \psi_2 - \partial_{xx} \psi_2 - \partial_{xx}\psi_1 \\
	\noalign{\smallskip}
-\psi_2
	\end{array}
	\right).
	$$

To build the differential operator $\mathcal{M}^\star$ (then we find $\mathcal{M}$), the goal is to apply some differential operator to the components of $\mathcal{L}^\star\Psi$ to obtain $\psi$.

	Let us introduce the operator $\mathcal{M}^\star_1:=\left(0 \quad - 1 \quad \partial_t+\partial_{xx}\right)$. Then, we obtain
	$$
\mathcal{M}^\star_1\circ\mathcal{L}^\star\psi=\partial_{xx}\psi_1 .
	$$
Taking $\mathcal{M}^\star_2:= - \mathcal{M}^\star_1- (1 \quad 0 \quad 0)$, we deduce
	$$
\mathcal{M}^\star_2\circ\mathcal{L}^\star\psi = \partial_t \psi_1 - q (\cdot) \psi_1 .
	$$
Using \eqref{eq:partialxq}, we define $\mathcal{M}^\star_3:=(- 2\partial_xq (\cdot ) )^{-1}[(-\partial_t + q (\cdot ) )\circ\mathcal{M}^\star_1+\partial_{xx}\circ\mathcal{M}^\star_2]$ in $\widehat{\omega}$. Hence, we have
	$$
\mathcal{M}^\star_3\circ\mathcal{L}^\star\psi = \partial_x \psi_1 + \frac{\partial_{xx}q(\cdot ) }{2\partial_x q(\cdot ) }\psi_1.
	$$
Setting $\mathcal{M}^\star_4:=\partial_{x}\circ\mathcal{M}^\star_2-\partial_t\circ\mathcal{M}^\star_3,$ it holds
	$$
\mathcal{M}^\star_4 \circ \mathcal{L}^\star \psi = - q(\cdot ) \partial_x \psi_1 - \partial_{x}q (\cdot ) \psi_1 - \frac{\partial_{xx}q (\cdot ) }{2\partial_x q(\cdot ) } \partial_t \psi_1 .
	$$
Again, using \eqref{eq:partialxq}, we consider $\mathcal{M}^\star_5:=- (q (\cdot ))^{-1}[\mathcal{M}^\star_4 + \frac{\partial_{xx}q(\cdot ) }{2\partial_xq(\cdot )}\mathcal{M}^\star_2] - \mathcal{M}^\star_3$. Then, we obtain
	$$
\mathcal{M}^\star_5 \circ \mathcal{L}^\star \psi = \frac{\partial_{x}q (\cdot )}{q (\cdot ) }\psi_1 .
	$$

Finally, if we take 
	$$
\mathcal{M}^\star:=\left(
	\begin{array}{c}
\displaystyle \frac{q (\cdot )}{\partial_{x}q(\cdot )}\mathcal{M}^\star_5 \\
	\noalign{\smallskip}
\mathcal{M}^\star_6
	\end{array}
	\right),
	$$
where $\mathcal{M}^\star_6:=(0 \quad 0 \quad -1)$, we obtain \eqref{eq:L*M*}. Therefore, we have also solved the algebraic problem~\eqref{eq:alg prob}. This ends the proof of Theorem \ref{th:fic}.
\end{proof}


\begin{remark}
The idea of algebraic resolvability for differential operators can be found in \cite[Section $2.3.8$]{gromov}. In the context of control theory, 
the Gromov algebraic resolvability was widely used, for instance in \cite{coron-lissy} for a Navier Stokes control system, in \cite{CORON_OLIVE} for first order quasi-linear hyperbolic systems and in \cite{Duprezlissy1,Duprezlissy2} for zero and first order coupled linear parabolic systems with a reduced number of distributed controls.
\end{remark}

\paragraph{Acknowledgements.}
\MD{The} second and third authors have been partially supported by  Grant~PID$2020$--$114976$GB--I$00$, funded by MCIN/AEI/$10.13039/501100011033$. \MD{The} third author was partially supported by  Grant IJC$2018$--$037863$-I funded by MCIN/AEI/$10.13039/501100011033$.

\end{document}